\documentclass[moor,nonblindrev]{informs1}

\usepackage{physics}
\usepackage{mathtools}
\usepackage{tensor}
\usepackage{comment}

\usepackage[mathscr]{euscript}   

\newcommand{\half}{\tfrac{1}{2}}

\newcommand{\ifft}{\mathop{\sf IFFT} \nolimits}

\newmuskip\pFqmuskip
\newcommand*\pFq[6][8]{%
  \begingroup 
  \pFqmuskip=#1mu\relax
  \mathcode`\,=\string"8000
  \begingroup\lccode`\~=`\
  \lowercase{\endgroup\let~}\pFqcomma
  {}_{#2}F_{#3}{\left[\genfrac..{0pt}{}{#4}{#5};#6\right]}%
  \endgroup
}
\newcommand{\pFqcomma}{\mskip\pFqmuskip}

\def\wscl{0.8}
\def\hscl{0.56}

\def\wsclx{0.8}
\def\hsclx{0.48}

\usepackage{natbib}
\NatBibNumeric
\def\bibfont{\small}%
\bibpunct[, ]{[}{]}{,}{n}{}{,}%

\definecolor{myblue}{RGB}{0, 20, 114}
\usepackage[colorlinks=true,breaklinks=true,bookmarks=true,urlcolor=myblue,
     citecolor=myblue,linkcolor=myblue,bookmarksopen=false,draft=false]{hyperref}

\def\EMAIL#1{\href{mailto:#1}{#1}}

\usepackage{algorithmic}
\usepackage{algorithm}

\TheoremsNumberedThrough     

\EquationsNumberedThrough    

\MANUSCRIPTNO{0} 

\def\theARTICLETOP{}
\def\theLRHSecondLine{\fs.7.9.\rm Article to be submitted to Annals of Operations Research}
\def\theRRHSecondLine{\fs.7.9.\rm Article to be submitted to Annals of Operations Research}

\begin{document}


\RUNAUTHOR{Zuk and Kirszenblat}

\RUNTITLE{Non-Preemptive Multi-Level Priority Queue}

\TITLE{Joint Queue-Length Distribution for the Non-Preemptive Multi-Server Multi-Level Markovian Priority Queue}

\ARTICLEAUTHORS{%
\AUTHOR{Josef Zuk}
\AFF{Defence Science and Technology Group, Melbourne, Australia,
\EMAIL{josef.zuk@defence.gov.au}}
\AUTHOR{David Kirszenblat}
\AFF{Defence Science and Technology Group, Melbourne, Australia,
\EMAIL{david.kirszenblat@defence.gov.au}}
} 

\ABSTRACT{%
Explicit results are obtained using simple and exact methods for the joint
queue-length distribution of the M/M/$c$ queue with an arbitrary number of
non-preemptive priority levels.
This work is the first to provide explicit results for the joint probability
generating function and joint probability mass function
for a general number of priority levels.
A fixed-point iteration is developed for the stationary balance equations,
which enables direct computation of the joint queue-length distribution.
A multi-variate probability generating function is also derived, from which the joint
probability mass function can be computed by means of a multi-dimensional
fast Fourier transform method.
}%

\MSCCLASS{Primary: 90B22; secondary: 60K25, 60J74}
\ORMSCLASS{Primary: Queues: Priority; secondary: Queues: Markovian}

\HISTORY{Date created: July 27, 2023. Last update: October 24, 2023.}

\maketitle

%

\section{Introduction}
\label{intro}
This work is concerned with the development of practical algorithms for computing
the joint queue-length distribution for the non-preemptive Markovian priority queue
with a general number of priority levels.
In the most recent edition of their textbook,
\citet{NP:Shortle18} remark that
`{\em the determination of stationary probabilities
in a non-preemptive Markovian system is an exceedingly difficult matter, well near
impossible when the number of priorities exceeds two}'.
\citet{NP:Elmelegy10} has also commented that `{\em given the immense
literature studying non-preemptive priority queuing systems, it is hard to find a simple
and exact method that calculates the performance measures of non-preemptive priority
systems with more than two priority levels}'.
The present discussion serves to fill this knowledge gap.

Previous work on the non-preemptive priority M/M/$c$ queue has been
reviewed recently in \citep{NP:Zuk23}.
The vast majority of effort concerning the joint distribution has focused on the two-level problem.
The single source of previous work addressing the joint queue-length distribution
for more than two priority levels comprises the papers of \citet{NP:Wignall73A} and \citet{NP:Wignall73}.
These consider single-server systems involving multiple queues
with distinct arrival rates ranked by priority level,
and with feedback -- requiring deterministic or probabilistic transitions between
queues before system exit is achieved.
The present problem can, in principle, be constructed as a special case of this scheme.
In \citep{NP:Wignall73A,NP:Wignall73}, equations that must be solved recursively
are presented  for the multivariate probability generating function (PGF); but
no general solution is provided, and manual solution becomes increasingly cumbersome
as the number of priority levels grows beyond a small number.
No actual probability mass functions (PMFs) are computed.
By contrast, we present an explicit closed-form expression for the joint PGF
given any number of priority levels, that lends itself to practical
numerical evaluation of the joint PMF.
In \citep{NP:Wignall73A}, explicit results are confined to the probability that
the system is empty, and the probability that a given queue is being served at a
random observation.
In \citep{NP:Wignall73}, the distribution of the maximum queue length during a busy
period in the presence of probabilistic feedback is also given.
The approach that we have adopted in the present paper may be viewed
as an extension of the method employed by \citet{NP:Cohen56} for the two-level case.
With it, we are able, in the words of \citet{NP:Neuts84}, `{\em to obtain actual results,
that is numbers and insight from numbers}'.

Priority-level numbers significantly greater than two are encountered in numerous real-world
applications, such as health care \citep{NP:Elalouf22}.
In a hospital emergency department (ED), arrivals are prioritized according to
patient acuity level \citep{NP:Hou20}. Most hospitals operate with at least five acuity levels.
In a forthcoming paper, the results obtained here will be applied to the ambulance ramping
problem \citep{NP:Almehdawe13},
in which arrivals to the ED by ambulance or as walk-ins are categorized into
three priority levels, corresponding to high, intermediate or low patient acuity.
The model is further complicated by the fact that
there are two arrival classes (ambulance and walk-in)
each with their own arrival rate and
each of which contain patients of multiple priority levels.
Thus, there is an entanglement between arrival classes and priority levels.

The number of servers is denoted by $c$ and the number of priority levels by $K$.
Each priority level is associated with a Poisson arrival rate $\lambda_\kappa$,
\mbox{$\kappa = 1,2,\ldots,K$},
leading to a total arrival rate
\mbox{$\lambda = \sum_{\kappa=1}^K\lambda_\kappa$}.
A common service rate $\mu$, associated with an exponential distribution,
is assumed for all priority levels.
Thus, the total traffic intensity is given by
\mbox{$r = \lambda/(c\mu)$}.
As we are interested in the state-state queue-length distribution,
possible values of $r$ will be limited to the ergodic region
\mbox{$r < 1$}.
The level traffic intensity for priority
\mbox{$\kappa = 1,2,\ldots,K$}
is defined as
\mbox{$r_\kappa = \lambda_\kappa/(c\mu)$},
so that
\mbox{$r = \sum_{\kappa=1}^K r_\kappa$}.
It is also convenient to introduce priority-level fractions
\mbox{$0\le \nu_\kappa \le 1$},
summing to unity, such that
\mbox{$r_\kappa = r\nu_\kappa$}.
Alternatively, the
\mbox{$\nu_\kappa \geq 0$}
may be chosen without constraint, provided we set
\mbox{$r_\kappa = r\nu_\kappa/\|\boldsymbol{\nu}\|_1$}.
This facilitates exploration of different distributions among the priority levels
for a given constant total load on the system.
To test the numerical performance of the algorithms developed here,
we shall fix $r$ and compute distributions for a random sample of vectors
\mbox{$\boldsymbol{\nu} = (\nu_1,\nu_2,\ldots,\nu_K)$}
chosen independently from the unit probability simplex in $K$-dimensions
\mbox{$\Delta_K \equiv \{\boldsymbol{\nu}\in \mathbb{R}^K: \|\boldsymbol{\nu}\|_1 = 1,
     \boldsymbol{\nu} \geq 0\}$}.
We also define the partial (per server) traffic intensity as
\mbox{$\rho \equiv \lambda/\mu$}
so that
\mbox{$r = c\rho$},
in line with the notation of \citep{NP:Gnedenko89}.

The rest of the paper is organized as follows: In Section~\ref{StatBalance}, we
set up the stationary balance equations for the model.
It is shown that, for any number of priority levels, they can be stated
in compact form comprising just a single equation.
This leads to a nearest-neighbour relationship among the elements of the joint PMF,
that is subsequently used as a diagnostic test of the computational procedures
developed later.
It is also shown that the balance equation can be solved directly for the joint PMF
by means of a fixed-point iteration (FPI).
While the FPI does not constitute an efficient computational algorithm, it
provides an important benchmark for verifying the correctness of the
vastly more efficient methodology discussed in the remainder of the paper.
In Section~\ref{PGF}, the general balance equation is used to derive
an explicit closed-form expression for the multi-variate
PGF of the joint distribution.
It is shown that the correct marginals follow from the joint PGF.
Section~\ref{FFT} describes a multi-dimensional
fast Fourier transform-based (FFT) method that
computes the joint PMF from the multi-variate PGF.
It is equally applicable to the marginal distributions.
Measures of performance for various diagnostic tests are presented
in Section~\ref{Tests}, the results of which are discussed there.
Conclusions follow in Section~\ref{Concl}.
Various technical details appear in the Appendices.

\section{Stationary Balance Equations}
\label{StatBalance}
\begin{figure}
\FIGURE
{\includegraphics[width=\wsclx\linewidth, height=\hsclx\linewidth]{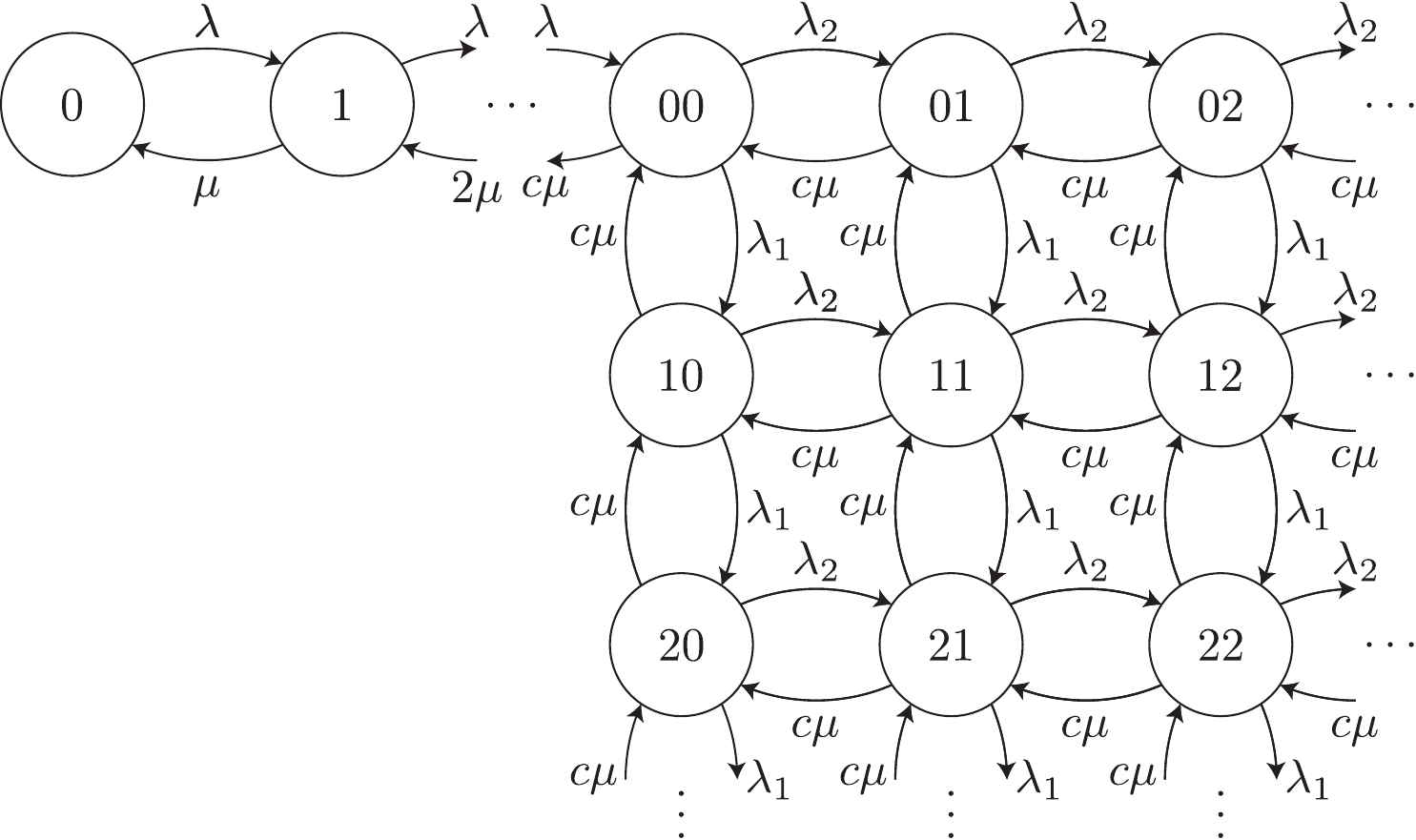}}
{\hphantom{x}\label{Markov}}
{Markov chain transitions for the two-level non-preemptive priority queue.}
\end{figure}

For ease of illustration, we shall begin by considering the concrete example of a system
with three priority levels: high, intermediate, low.
Let the single-subscripted quantity $p_k$ denote the probability that the queue is empty and that
$k$ of the $c$ servers are busy.
Let the triple-subscripted quantity $p_{\ell mn}$ denote the probability that there are $\ell$ high-priority
clients, $m$ intermediate-priority clients and $n$ low-priority client in the queue, and that all servers are busy.
Since, apart from the identification
\mbox{$p_c = p_{000}$},
the probabilities $p_k$ and $p_{\ell mn}$ are exhaustive and mutually exclusive, we must have that
\begin{equation}
\sum_{k = 0}^{c-1} p_k + \sum_{\ell,m,n = 0}^\infty p_{\ell mn} = 1 \;.
\label{Norm}
\end{equation}
Figure~\ref{Markov} displays the state transition diagram for the two-level case, which is easily drawn
and captures all the qualitative features of the general problem.
With the definitions above, the stationary balance equations for the three-level system are given by
%
%
\begin{align}
\begin{aligned}
&(a)\quad &\lambda p_k                 &= (k+1)\mu p_{k+1} &\text{for} \quad &0 \leq k \leq c \;, \\
&(b)\quad &(\lambda + c\mu)p_{000}     &= \lambda p_{c-1} + c\mu(p_{100} + p_{010} + p_{001}) & &\\
&(c)\quad &(\lambda + c\mu)p_{00n}     &= \lambda_3 p_{00n-1} + c\mu(p_{10n} + p_{01n} + p_{00n+1}) &\text{for} \quad &n \geq 1 \;, \\
&(d)\quad &(\lambda + c\mu)p_{0mn}     &= \lambda_2 p_{0m-1n} + \lambda_3 p_{0mn-1} +c\mu(p_{1mn} + p_{0m+1n}) &\text{for} \quad &m \geq 1, n \geq 0 \;, \\
&(e)\quad &(\lambda + c\mu)p_{\ell mn} &= \lambda_1 p_{\ell-1mn} + \lambda_2 p_{\ell m-1n} + \lambda_3 p_{\ell mn-1} +  c\mu p_{\ell+1mn} &\text{for} \quad &\ell \ge 1, m,n \geq 0 \;.
\end{aligned}
\label{Balance}
\end{align}
We adopt the convention that negative valued subscripts correspond to zero probability.
It follows immediately from (a) that
\begin{equation}
p_k = \frac{1}{k!}\left(\frac{\lambda}{\mu}\right)^k p_0 \;,
\end{equation}
for
\mbox{$k = 0,1,\ldots,c$}.
Therefore,
\begin{equation}
p_{000} \equiv p_c = \frac{1}{c!}\left(\frac{\lambda}{\mu}\right)^c p_0 \;.
\end{equation}
With the observation that
\mbox{$\lambda p_{c-1} = c\mu p_c = c\mu p_{000}$},
we can recast the collection of balance equations above for the $p_{\ell mn}$ into the single equation
\begin{equation}
\begin{split}
(\lambda + c\mu)p_{\ell mn} &= \lambda_1 p_{\ell-1mn} + \lambda_2 p_{\ell m-1n} + \lambda_3 p_{\ell mn-1} \\
&    {}+ c\mu\left(p_{\ell+1mn} + \delta_{0\ell}p_{\ell m+1n} + \delta_{0\ell}\delta_{0m}p_{\ell mn+1} +
     \delta_{0\ell}\delta_{0m}\delta_{0n}p_{\ell mn}\right) \;.
\end{split}
\label{Single}
\end{equation}
One should note that (\ref{Single}) is homogeneous in $p_{\ell mn}$,
which means that its solution is determined only up to an arbitrary scale factor.
It is convenient to set
\mbox{$\tilde{p}_{\ell mn} \equiv p_{\ell mn}/p_{000}$}
so that
\mbox{$\tilde{p}_{000} = 1$}.
The normalization condition (\ref{Norm}) then implies that
\begin{equation}
\frac{1}{p_{000}} = \sum_{k = 0}^{c-1}\frac{c!}{k!}\left(\frac{\mu}{\lambda}\right)^{c-k} +
     \sum_{\ell,m,n=0}^\infty \tilde{p}_{\ell mn} \;.
\label{pzero}
\end{equation}
Thus, we may solve (\ref{Single}) for $\tilde{p}_{\ell mn}$ by setting
\mbox{$\tilde{p}_{000} = 1$},
and subsequently determine $p_{000}$ from (\ref{pzero}),
so that the probabilities can be recovered via
\mbox{$p_{\ell mn} = p_{000}{\cdot}\tilde{p}_{\ell mn}$}.

Next, let us write
\begin{equation}
\Phi(\ell) \equiv \sum_{m,n}^\infty p_{\ell mn} \;,
\end{equation}
and observe that we also have
\begin{equation}
\sum_{m,n=0}^\infty p_{\ell m-1n} = \sum_{m,n=0}^\infty p_{\ell mn-1} = \Phi(\ell) \;.
\end{equation}
Moreover,
\begin{align}
\begin{aligned}
\sum_{m,n=0}^\infty \delta_{0\ell}p_{\ell m+1n}                       &= \delta_{0\ell}\Bigl[\Phi(0) - \sum_{n=0}^\infty p_{00n}\Bigr] \;, \\
\sum_{m,n=0}^\infty \delta_{0\ell}\delta_{0m}p_{\ell mn+1}            &= \delta_{0\ell}\Bigl[\sum_{n=0}^\infty p_{00n} - p_{000}\Bigr] \;, \\
\sum_{m,n=0}^\infty \delta_{0\ell}\delta_{0m}\delta_{0n}p_{\ell mn} &= \delta_{0\ell}p_{000}\;,
\end{aligned}
\end{align}
and so we see that these three contributions sum to $\Phi(0)$.
Accordingly, (\ref{Single}) implies that
\begin{equation}
(\lambda_1 + c\mu)\Phi(\ell) = \lambda_1\Phi(\ell-1) + c\mu\left[\Phi(\ell+1) + \delta_{0\ell}\Phi(0)\right] \;,
\label{Phi}
\end{equation}
for
\mbox{$\ell = 0,1,\ldots$}.
If we let
\mbox{$\Delta\Phi(\ell) \equiv \Phi(\ell) - \Phi(\ell-1)$},
then
\mbox{$\Delta\Phi(\ell+1) = r_1\Delta\Phi(\ell)$}
for
\mbox{$\ell = 1,2,\ldots$},
subject to
\mbox{$\Delta\Phi(1) = -(1-r_1)\Delta\Phi(0)$}
and where
\mbox{$r_1 = \lambda_1/(c\mu)$}.
This is solved by
\mbox{$\Delta\Phi(\ell) = \varphi(\ell){\cdot}\Delta\Phi(0)$},
with
\mbox{$\varphi(\ell) \equiv -(1 - r_1)r_1^{\ell-1}$}.
It is equivalent to the first-order recurrence equation
\begin{equation}
\Phi(\ell) = \Phi(\ell-1) + \varphi(\ell){\cdot}\Phi(0) \;,
\end{equation}
for
\mbox{$\ell = 1,2,\ldots$}
which, in turn,  is solved by
\begin{equation}
\Phi(\ell) = \Phi(0){\cdot}\Bigl[1 + \sum_{k= 1}^\ell\varphi(k)\Bigr] = \Phi(0)r_1^\ell \;,
\label{PhiSolve}
\end{equation}
for
\mbox{$\ell = 0,1,\ldots$}.
The quantity $\Phi(\ell)$ represents the probability that there are $\ell$ high-priority
clients in the queue and that all servers are busy.
For
\mbox{$\ell \geq 1$},
it simply gives the probability of $\ell$ high-priority clients in the queue without qualification
because having one or more high-priority clients in the queue implies that all servers must be busy.
Let $\Psi(\ell)$ denote the probability of $\ell$ high-priority clients in the queue.
Then, for
\mbox{$\ell \geq 1$},
\mbox{$\Psi(\ell) = \Phi(\ell) = \Phi(0)r_1^\ell$}
and
\begin{equation}
1 = \sum_{\ell=0}^\infty \Psi(\ell) = \Psi(0) + \sum_{\ell=1}^\infty \Phi(\ell) \;,
\end{equation}
from which it follows that
\mbox{$\Psi(0) = 1 - [r_1/(1 - r_1)]{\cdot}\Phi(0)$}.
Thus, we may write
\begin{equation}
\Psi(\ell) =  \left[1 - \Phi(0)/(1-r_1)\right]\delta_{\ell 0} + \Phi(0)r_1^\ell \;.
\label{PsiSolve}
\end{equation}
The normalization condition (\ref{Norm}) implies that
\begin{equation}
\sum_{\ell=0}^\infty \Phi(\ell) = 1 - P_{\text{NW}}\;, \quad P_{\text{NW}} \equiv \sum_{k=0}^{c-1} p_k \;,
\end{equation}
where $P_{\text{NW}}$ represents the no-wait probability.
We may equate this with the summation over $\Phi(\ell)$ as given by (\ref{PhiSolve}) to obtain
\mbox{$\Phi(0) = (1-r_1)(1 - P_{\text{NW}})$}.
Hence, (\ref{PsiSolve}) can be re-expressed in the more direct form
\begin{equation}
\Psi(\ell) =  P_{\text{NW}}\delta_{\ell 0} + (1 - P_{\text{NW}})(1-r_1)r_1^\ell \;.
\end{equation}
This result applies equally to the general problem with an arbitrary number of priority levels.
The only undetermined component in the foregoing discussion is the value of $p_0$, which represents
the probability that the system is empty, and is given by \citep{NP:Gnedenko89}
\begin{equation}
\frac{1}{p_0} = \sum_{k=0}^c \frac{\rho^k}{k!} + \frac{\rho^{c+1}}{c!(c - \rho)} \;,
\end{equation}
where
\mbox{$\rho \equiv \lambda/\mu$}.
Knowledge of $p_0$ allows us to obtain
\mbox{$p_{000} = p_c$}.
Another way of characterizing it is as follows:
Let the random variable $\mathcal{N}_{\text{sys}}$ represent the number of clients in the system.
It is well-known that, for some constant $A$,
\mbox{$\Pr(\mathcal{N}_{\text{sys}} = n) = Ar^n$}
for all
\mbox{$n \geq c$}.
By construction,
\mbox{$p_{000} = \Pr(\mathcal{N}_{\text{sys}} = c)$}.
We also have
\mbox{$P_{\text{NW}} = \Pr(\mathcal{N}_{\text{sys}} \leq c-1)$}.
Hence,
\begin{align}
\begin{aligned}
P_{\text{NW}} &= \sum_{n=0}^{c-1}\Pr(\mathcal{N}_{\text{sys}} = n)\!\!\!
     &= 1 &- \sum_{n=c}^\infty\Pr(\mathcal{N}_{\text{sys}} = n) \\
     &= 1 - \frac{Ar^c}{1-r}
     &= 1 &- \frac{1}{1-r}\Pr(\mathcal{N}_{\text{sys}} = c) \;,
\end{aligned}
\end{align}
which leads to the result
\mbox{$p_{000} = (1-r)(1-P_{\text{NW}})$}.

The stationary balance equations as given in (\ref{Single}) are easily extended to the general problem of
an arbitrary number $K$ of priority levels.
In order to establish a compact expression, we introduce the lattice vectors
\mbox{$\mathbf{n} \equiv [n_1,n_2,\ldots,n_K] \in \mathbb{Z}^K$},
and let $\mathbf{e}_\kappa$,
\mbox{$\kappa = 1,2,\ldots,K$}
denote the standard unit Cartesian coordinate basis vectors in $\mathbb{Z}^K$.
Then the stationary balance equations for $K$ priority levels can be expressed as
\begin{equation}
(1 + r)p_\mathbf{n} = \prod_{j=1}^K\delta_{0n_j}{\cdot}p_{\mathbf{n}}
     + \sum_{\kappa=1}^K \biggl[ r_\kappa p_{\mathbf{n} - \mathbf{e}_\kappa}
     + \prod_{j=1}^{\kappa-1}\delta_{0n_j}{\cdot}p_{\mathbf{n} + \mathbf{e}_\kappa}\biggr] \;,
\label{Vector}
\end{equation}
where we recall the convention that
\mbox{$p_\mathbf{n} \equiv 0$}
if
\mbox{$n_\kappa < 0$}
for any
\mbox{$\kappa = 1,2,\ldots,K$}.
If we introduce the boundary set\footnote{Strictly speaking, this is the boundary and
     beyond or, equivalently, the set of non-interior points.}
\begin{equation}
\mathcal{B} \equiv \{\mathbf{n} \in \mathbb{Z}^K: n_\kappa \leq 0 \text{ for some } \kappa = 1,2,\ldots,K\} \;,
\end{equation}
then, for all
\mbox{$\mathbf{n} \notin \mathcal{B}$},
we have the interior stationary balance equations
\begin{equation}
p_\mathbf{n} = \frac{1}{1+r}\biggl[p_{\mathbf{n} + \mathbf{e}_1} +
     \sum_{\kappa=1}^K r_\kappa p_{\mathbf{n} - \mathbf{e}_\kappa}\biggr] \;.
\label{Interior}
\end{equation}
Thus every probability in the interior region is a positive weighted sum of its lower nearest neighbours
plus its upper highest priority neighbour.
While it is not a numerically stable proposition to attempt to first solve for $p_\mathbf{n}$ on the
boundary set $\mathcal{B}$ and then use (\ref{Interior}) to propagate the solution into the interior region,
the relationship (\ref{Interior}) does provide a robust diagnostic test of where a candidate solution for
$p_\mathbf{n}$ behaves as it should.
We shall call this the nearest-neighbour test.
The two-dimensional instance of this test was used extensively for the two-level problem in \citep{NP:Zuk23}.
The wait-conditional joint queue-length PMF is given by
\mbox{$P(\mathbf{n}) = (1-r)\tilde{p}_\mathbf{n}$}.
The full unconditional joint PMF is then expressed as
\begin{align}
\begin{aligned}
P_{\text{full}}(\mathbf{n}) &= P_{\text{NW}}\delta(\mathbf{n}) + (1-P_{\text{NW}})P(\mathbf{n}) \\
     &= P_{\text{NW}}\delta(\mathbf{n}) + (1-P_{\text{NW}})(1-r)\tilde{p}_\mathbf{n} \;,
\end{aligned}
\label{Pfull}
\end{align}
where
\mbox{$\delta(\mathbf{n}) \equiv \prod_{\kappa=1}^K\delta_{0n_\kappa}$},
and recalling that
\mbox{$\tilde{p}_\mathbf{0} = 1$}.

The system (\ref{Vector}) can be solved directly using an FPI.
However, this necessitates truncation of the problem to a maximum queue size for each priority level.
Consider the space $\mathcal{A}_\infty^K$ of countably infinite matrices in $K$-dimensions, whose
elements are non-negative and indexed by
the $K$-dimensional vector
\mbox{$\mathbf{n} = [n_1,n_2,\ldots,n_K]$},
with
\mbox{$n_\kappa = 0,1,\ldots$},
for each
\mbox{$\kappa = 1,2\ldots,K$}.
In other words,
\mbox{$\mathbf{n}\in {\mathbb{N}_0^K}$}
--- the space of $K$-tuples of non-negative integers.
We define the matrix mapping
\mbox{$\mathcal{M}: \mathcal{A}_\infty^K \to \mathcal{A}_\infty^K$}
by
\begin{equation}
A' = \mathcal{M}(A)\;, \quad A'_\mathbf{n} = \frac{1}{1+r}\left[
     \prod_{j=1}^K\delta_{0n_j}{\cdot}A_{\mathbf{n}}
     + \sum_{\kappa=1}^K \biggl( r_\kappa A_{\mathbf{n} - \mathbf{e}_\kappa}
     + \prod_{j=1}^{\kappa-1}\delta_{0n_j}{\cdot}A_{\mathbf{n} + \mathbf{e}_\kappa}\biggr)\right] \;.
\end{equation}
We observe that the mapping $\mathcal{M}$ is sum preserving:
\begin{equation}
A' = \mathcal{M}(A) \Rightarrow \sum_{\mathbf{n}\in{\mathbb{N}_0^K}} A'_\mathbf{n}
     = \sum_{\mathbf{n}\in{\mathbb{N}_0^K}} A_\mathbf{n} \;,
\end{equation}
assuming that the sum is finite.
This is equivalent to
\mbox{$\|A'\|_1 = \|A\|_1$}
for the $L_1$ matrix norm,
and follows from (\ref{Phi}), which trivially extends to the general multi-level case.

We shall consider the collection of probabilities $p_\mathbf{n}$ as defining a matrix
\mbox{$\mathsf{P} \in \mathcal{A}_\infty^K$}.
The $L_1$ and $L_\infty$ norms are given, respectively, by
\begin{equation}
\|\mathsf{P}\|_1 \equiv \sum_{\mathbf{n}} |\mathsf{P}_\mathbf{n}| \;, \quad
\|\mathsf{P}\|_\infty \equiv \max_{\mathbf{n}} |\mathsf{P}_\mathbf{n}| \;.
\end{equation}
Both these norms exist and are finite for our matrix of
(non-negative) probabilities $\mathsf{P}$.

To solve the FPI directly, we must truncate the queue length in each dimension
({\it i.e.}\ priority level) to a finite maximum size,
\mbox{$n_\kappa \leq N_\kappa^{\text{max}}$}
for each
\mbox{$\kappa = 1,2,\ldots,K$}.
For the sake of simplicity, we shall truncate in each of the $K$ dimensions
to a common value $N_{\text{max}}$.
For the FPI, the matrix $\mathsf{P}$ will be initialized to zero, except for
\mbox{$\mathsf{P}_\mathbf{0} = 1$},
on the grid
\mbox{$[-1,N_{\text{max}}+1]^K$},
but subsequently updated only on the sub-grid
\mbox{$[0,N_{\text{max}}]^K$}.
Finite-size truncation will cause inevitable leakage of probability at each
step of the FPI. However, we can use the $L_1$ invariance to quantify the
magnitude of the probability leakage, and adjust for it by amortizing uniformly
over all matrix elements.
Also, after each FPI step, we scale the matrix $\mathsf{P}$ to preserve
\mbox{$\mathsf{P}_\mathbf{0} = 1$}.
The correct normalization is applied at the end.
The details are provided in Algorithm~\ref{FixedPt}.

\begin{algorithm}
\caption{Fixed-point iteration.}
\label{FixedPt}
\begin{algorithmic}[1]
\REQUIRE $(r_1,\ldots,r_K), N_{\text{max}}$
\ENSURE $\mathsf{P}$
\STATE
\COMMENT{Parameters:}
\STATE
$\epsilon_{\text{tol}} = 10^{-9}$
\STATE
\COMMENT{Initialization:}
\STATE
$\Delta = \infty$
\STATE
$\mathsf{P}_\mathbf{n} \leftarrow \delta(\mathbf{n})$
\WHILE{$\Delta > \epsilon_{\text{tol}}$}
  \STATE \ \ \
  \COMMENT{Iteration:}
  \STATE \ \ \
  $\mathsf{P}' = \mathcal{M}(\mathsf{P})$
  \STATE \ \ \
  \COMMENT{Probability leakage amortization:}
  \STATE \ \ \
  $p_{\text{leak}} = \|\mathsf{P}\|_1 - \|\mathsf{P'}\|_1 > 0$
  \STATE \ \ \
  $\mathsf{P}'     \leftarrow \mathsf{P}' + p_{\text{leak}}/(N_{\text{max}}+1)^K$
  \STATE \ \ \
  \COMMENT{Renormalization:}
  \STATE \ \ \
  $\mathsf{P}' = \mathsf{P}'/\mathsf{P}'_\mathbf{0}$
  \STATE \ \ \
  \COMMENT{Convergence:}
  \STATE \ \ \
  $\Delta = \|\mathsf{P}' - \mathsf{P}\|_\infty$
  \STATE \ \ \
  $\mathsf{P} \leftarrow \mathsf{P}'$
\ENDWHILE
\STATE
\COMMENT{Finalization:}
\STATE
$\mathsf{P} \leftarrow (1-r){\cdot}\mathsf{P}$
\end{algorithmic}
\end{algorithm}

Apart from the requirement that the queue length for each priority level be truncated
to a finite maximum value,
the FPI method has other disadvantages, as it suffers from the `curse of dimensionality',
and convergence becomes very slow as the total traffic intensity $r$ gets close to unity.
The slow convergence in this region is exacerbated by the fact that the maximum queue size needs
to be set quite large, as probabilities of large queue sizes become non-negligible.
Its role in the present discussion is that it serves as a benchmark for verifying the
correctness for the much more efficient FFT method developed in the following sections.
Not only will it confirm that the joint queue-length distribution
has been computed correctly, there will be an implied confirmation
of the multi-variate PGF on which the FFT method rests.
Results of this comparison are illustrated in Figure (\ref{TestRepFPIFFTJointLevelPMF})
and discussed later on in the results section.
It is difficult to use Monte Carlo (MC) simulation to check a multi-variate distribution
in its entirety due to size requirements and lack of appropriate hypothesis tests.
We have, however, performed the more straightforward task of
checking marginal distributions against discrete-event MC simulation.
The FPI method does have the singular advantage that one may expect it to be robust
to the incorporation into the model of a variety of complications
({\it e.g.}\ unequal service rates, or more complex queue disciplines),
whereas the closed-form PGF may not survive such modifications.

\section{Probability Generating Function}
\label{PGF}
Generalizing the approach of \citet{NP:Cohen56},
let us introduce the collection of multi-variate functions of $K-1$ continuous variables
\begin{equation}
G_\ell(\mathbf{u}) \equiv P_0{\cdot}\sum_{m_1,\ldots,m_{K-1}=0}^\infty
     \tilde{p}_{\ell m_1\cdots m_{K-1}}u_1^{m_1}\cdots u_{K-1}^{m_{K-1}} \;,
\end{equation}
for
\mbox{$\ell = 0,1,\ldots$},
where
\mbox{$\mathbf{u} \equiv [u_1,\ldots,u_{K-1}]$}.
With the identification
\mbox{$P_0 \equiv 1-r$},
we have
\begin{equation}
G_0(\mathbf{0}) = 1-r \;, \quad \sum_{\ell=0}^\infty G_\ell(\mathbf{1}) = 1 \;,
\end{equation}
and
\mbox{$G_\ell(\mathbf{u})$}
represents the PGF for the wait-conditional joint queue-length PMF,
such that
\begin{equation}
P(\mathbf{n}) = \prod_{\kappa=1}^{K-1}\frac{1}{n_{\kappa+1}!}
     \left.\frac{\partial^{n_{\kappa+1}}}{\partial u^{n_{\kappa+1}}_\kappa}{\cdot}
     G_{n_1}(u_1,\ldots,u_{K-1})\right|_{\mathbf{u}=\mathbf{0}} \;.
\end{equation}
Then, the PGF for the full unconditional PMF is given by
\begin{equation}
G_{\text{full},\ell}(\mathbf{u}) = P_{\text{NW}}\delta_{0\ell} + (1-P_{\text{NW}})G_\ell(\mathbf{u}) \;.
\end{equation}

Summing over the stationary balance equations yields the linear recurrence relations
\begin{equation}
G_{\ell+1} + [\alpha - (1+r)]G_\ell + r_1G_{\ell-1} = 0 \;,
\label{Grecur}
\end{equation}
for
\mbox{$\ell \geq 1$},
where
\begin{equation}
\alpha \equiv \mathbf{r}\cdot\mathbf{u} = \sum_{\kappa=1}^{K-1} r_{\kappa+1} u_\kappa \;.
\label{alpha}
\end{equation}
The characteristic equation reads
\begin{equation}
\lambda^2 + [\alpha - (1+r)]\lambda + r_1 = 0 \;,
\end{equation}
and is solved by
\mbox{$\lambda = \lambda_\pm$}
with
\begin{equation}
\lambda_\pm = \half\left[1 + r - \alpha \pm\sqrt{(1+r-\alpha)^2 - 4r_1}\right] \;.
\end{equation}
When
\mbox{$\mathbf{u} = \mathbf{1} \equiv [1,\ldots,1]$},
in which case
\mbox{$\alpha = r - r_1$},
we require, from the foregoing discussion of the high-priority marginal, that
\mbox{$G_\ell \propto r_1^\ell$},
and we have
\mbox{$\lambda_- = r_1$}
whereas
\mbox{$\lambda_+ = 1$}.
It follows immediately that
\begin{equation}
G_\ell(\mathbf{u}) = G_0(\mathbf{u})\lambda_-^\ell(\mathbf{u}) \;,
\label{G0Lambda}
\end{equation}
for
\mbox{$\ell = 0,1,\ldots$}.
It remains to solve the problem for
\mbox{$\ell = 0$}.

For this purpose, let us write
\mbox{$\mathbf{n} \equiv (\ell,\mathbf{m})$}
with
\mbox{$\mathbf{m} = (m_1,\ldots,m_{K-1}) \in \mathbb{N}_0^{K-1}$}
and
\mbox{$q_\mathbf{m}^\ell \equiv p_{\ell\mathbf{m}}$}.
Then, for
\mbox{$\ell = 0$},
and on setting
\begin{equation}
Q_\mathbf{m} \equiv \sum_{j=1}^{K-1}q^0_{\mathbf{m}+\mathbf{e}_j}{\cdot}
     \prod_{i=1}^{j-1}\delta_{0m_i} \;,
\end{equation}
the birth-death equations (\ref{Vector}) may be cast as
\begin{equation}
(1+r)q^0_\mathbf{m} = \prod_{j=1}^{K-1}\delta_{j0m_j}{\cdot}q^0_\mathbf{m}
     + \sum_{j=1}^{K-1}r_{j+1}q^0_{\mathbf{m} - \mathbf{e}_j} + q^1_\mathbf{m}
     + Q_\mathbf{m} \;.
\end{equation}
Summing over the multi-indices $\mathbf{m}$ with the corresponding
powers of $\mathbf{u}$, we see that
the stationary balance equations (\ref{Vector}) require that
\begin{equation}
(1+r)G_0(\mathbf{u}) = G_0(\mathbf{0}) + \alpha G_0(\mathbf{u})
     + G_1(\mathbf{u}) + \sum_{\mathbf{m}\in\mathbb{N}_0^{K-1}}
     Q_{\mathbf{m}}{\cdot}\prod_{j=1}^{K-1} u_j^{m_j} \;.
\end{equation}
%
Noting that
\mbox{$G_1(\mathbf{u}) =\lambda_- G_0(\mathbf{u})$},
and using the identity
\mbox{$\lambda_+ + \lambda_- = 1 + r -\alpha$},
we arrive at
\begin{equation}
\lambda_+ G_0(\mathbf{u}) = G_0(\mathbf{0}) + \sum_{\mathbf{m}\in\mathbb{N}_0^{K-1}}
    Q_{\mathbf{m}}{\cdot}\prod_{j=1}^{K-1} u_j^{m_j} \;,
\end{equation}
which may be further manipulated to yield
\begin{align}
\begin{aligned}
&\lambda_+ G_0(u_1,\ldots,u_{K-1}) \\
& = \sum_{k=1}^{K-1}\frac{1}{u_k}\left[G_0(0,\ldots,0,u_k,\ldots,u_{K-1}) -
     G_0(0,\ldots,0,u_{k+1},\ldots,u_{K-1})\right] + G_0(0,\ldots,0) \\
& = \frac{1}{u_1}G_0(u_1,\ldots,u_{K-1}) + \sum_{k=2}^{K-1} \left(\frac{1}{u_k} - \frac{1}{u_{k-1}}\right)
     G_0(0,\ldots,0,u_k,\ldots,u_{K-1}) + \left(1 - \frac{1}{u_{K-1}}\right)G_0(0,\ldots,0) \;.
\end{aligned}
\label{G0Req}
\end{align}
On formally setting
\mbox{$u_K \equiv 1$},
this may be simplified as
\begin{equation}
G_0(u_1,\ldots,u_{K-1}) = \frac{1}{1/u_1 - \lambda_+(u_1,\ldots,u_{K-1})}\sum_{k=2}^K
     \left(\frac{1}{u_{k-1}} - \frac{1}{u_k}\right)
     G_0(0,\ldots,0,u_k,\ldots,u_{K-1}) \;.
\label{G0U}
\end{equation}
It should be noted that
\mbox{$G'_0(u_k,\ldots,u_{K-1}) \equiv G_0(0,\ldots,0,u_k,\ldots,u_{K-1})$}
is the result for an appropriately aggregated $(K-k+1)$-dimensional problem,
obtained by treating the $k$ highest priority levels as a single high priority level
with level traffic intensity
\mbox{$r_{\text{hi}} = \sum_{\kappa=1}^k r_\kappa$}.
One may also note that
\mbox{$G'_0() = G_0(0,\ldots,0) = P_0$}.

One immediate consequence of (\ref{G0U}) is that the marginal distribution of the
aggregation of the top $p$ priority levels is geometric, for all
\mbox{$p = 1,2,\ldots,K$}.
From (\ref{G0Lambda}), the PGF for this distribution is
\begin{equation}
G_{\text{agg}}^{(p)}(u) \equiv \sum_{\ell=0}^\infty u^\ell G_\ell(\mathbf{u}^{(p)})
     = \frac{G_0(\mathbf{u}^{(p)})}{1 - u\lambda_-(\mathbf{u}^{(p)})} \;,
\label{Gaggp}
\end{equation}
where
\mbox{$\mathbf{u}^{(p)} \equiv (u{\cdot}\mathbf{1}_{p-1}, \mathbf{1}_{K-p})$}.
Only a single term survives the summation in (\ref{G0U}) in the evaluation of
\mbox{$G_0(\mathbf{u}^{(p)})$}, to produce
\begin{equation}
G_{\text{agg}}^{(p)}(u) = \frac{1-u}{[1 - u\lambda_+(\mathbf{u}^{(p)})]
     [1 - u\lambda_-(\mathbf{u}^{(p)})]}
     G_0(\mathbf{0}_{p-1},\mathbf{1}_{K-p}) \;.
\end{equation}
Now, we have
\begin{equation}
\alpha^{(p)} \equiv \mathbf{r}\cdot\mathbf{u}^{(p)} = (\sigma_p - \sigma_1)u
     + r - \sigma_p \;.
\end{equation}
where
\mbox{$\sigma_p \equiv \sum_{k=1}^p r_k$}.
It follows that we can write
\begin{equation}
\lambda_\pm(\mathbf{u}^{(p)}) = \half\left[1 + r_{\text{agg}} - r_{\text{lo}}u
     \pm\sqrt{(1 + r_{\text{agg}} - r_{\text{lo}}u)^2 - 4r_{\text{hi}}}\right] \;,
\end{equation}
with
\begin{equation}
r_{\text{lo}} = \sigma_p - \sigma_1 \;, \quad
     r_{\text{hi}} = \sigma_1 \;, \quad
     r_{\text{agg}} = r_{\text{lo}} + r_{\text{hi}} = \sigma_p \;.
\end{equation}
Hence,
\begin{equation}
\frac{1-u}{(1 - u\lambda_+)(1 - u\lambda_-)} = \frac{1}{1 - \sigma_p u} \;,
\end{equation}
and the desired result follows, consistent with intuitive expectations.

The solution of the system (\ref{G0U}) is derived in Appendix~\ref{AppPGF}.
Here, we explain the overall strategy:
We begin by noting that the recurrence relation (\ref{Grecur}) is similar to
equation (1.2.10) of \citet{NP:Cohen56},
but with his analogue of $r_2 u_1$ replaced by our $\alpha$ defined in (\ref{alpha}).
The analogue of (1.2.12) still holds, but the analogue of (1.2.11)
acquires additional complicating terms of the form
$G_0(0,\ldots,0,u,v,w,\ldots)$.
However, these can be dealt with by means of a simple aggregation trick.
This is because the leading zeros mean that we looking at probabilities of
various numbers of low priority entities in the queue and no entities with
the highest $n$ priorities.
But this is the same as no entities of a single highest priority level
aggregated from the highest $n$ priorities.
Thus, we have to solve a lower-dimensional problem.
The relevant observation can be expressed, for example in a 4-level problem, as
\mbox{$G_0(0, v, w) = G'_0(v, w)$},
\mbox{$G_0(0, 0, w) = G'_0(w)$},
and of course
\mbox{$G_0(0, 0, 0) = G'_0() = P_0$}.
The PGFs with lesser numbers of arguments refer to appropriately aggregated
lower-dimensional problems. Therefore, one has to start with Cohen's
two-dimensional problem and successively work upwards to the full
$K$-dimensional problem.
The recurrence relations take on a simpler, more transparent structure,
that is amenable to explicit solution, when one changes notation
by writing the arguments of the PGFs backwards.
Accordingly, it is convenient to present the result for $G_0(.)$ in terms of
new variables
\mbox{$(z_1,z_2,\ldots,z_{K-1})$}
where ascending indices are associated with ascending priority level.
Thus
\mbox{$z_k = u_{K-k}$},
in which case
\begin{equation}
G_0(0,\ldots,0,u_k,\ldots,u_{K-1}) = G_0(0,\ldots,0,z_{K-k},\ldots,z_2,z_1)
\end{equation}
or, equivalently,
\mbox{$G'_0(u_k,\ldots,u_{K-1}) = G'_0(z_{K-k},\ldots,z_2,z_1)$}.
The result is then given for a new version of $G_0(.)$
re-parameterized according to
\mbox{$G_0(z_1,\ldots,z_{K-1}) \equiv G'_0(z_{K-1},\ldots,z_2,z_1)$},
which is just $G'_0(.)$ with its arguments flipped.

A self-contained summary of Appendix~\ref{AppPGF} is as follows:
To construct the PGF for the problem of $K$ non-preemptive priority levels,
we consider the joint PMF $p_{\ell mn\cdots}$,
where $\ell$ enumerates the number in the queue of the highest priority level,
and write the PGF as
\begin{equation}
G_\ell(z_1,\ldots,z_{K-1}) \equiv \sum_{m_1,\ldots,m_{K-1} = 0}^\infty p_{\ell m_1\cdots m_{K-1}}
     z_{K-1}^{m_1}\cdots z_1^{m_{K-1}} \;,
\end{equation}
so that $z_1$ is associated with the {\em lowest} priority level and $z_{K-1}$ with the
{\em next-to-highest} priority level.
We focus our attention on the
\mbox{$\ell = 0$}
component
\mbox{$G_0(z_1,\ldots,z_{K-1})$}, and write
\mbox{$P_0 \equiv G_0(0,\ldots,0)$}.
Then, we have
\begin{equation}
G_0(z_1,\ldots,z_{K-1}) = P_0{\cdot}\prod_{\kappa=1}^{K-1}\frac{1 - z_\kappa\zeta_+(z_1,\ldots,z_{\kappa-1})}
     {1 - z_\kappa\zeta_+(z_1,\ldots,z_{\kappa})} \;.
\label{G0}
\end{equation}
The
\mbox{$\zeta_\pm(z_1,\ldots,z_\kappa)$}
solve the quadratic equation
\begin{equation}
\zeta^2 + [\beta(z_1,\ldots,z_\kappa) - (1+r)]\zeta + \sigma_{K-\kappa} = 0 \;,
\end{equation}
for
\mbox{$\kappa = 1,2,\ldots,K-1$}.
We have defined
\mbox{$\sigma_\kappa \equiv \sum_{k = 1}^\kappa r_k$},
so that
\mbox{$\sigma_{K-\kappa} = \sum_{k = 1}^{K-\kappa} r_k$}
and
\mbox{$\sigma_K = r$}.
Also,
\begin{equation}
\beta(z_1,\ldots,z_\kappa) \equiv \sum_{k = 1}^\kappa z_k r_{K+1-k} \;,
\end{equation}
for
\mbox{$\kappa = 1,2,\ldots,K-1$}.
We adopt the conventions that
\mbox{$z_0 \equiv 1$},
\mbox{$\beta() \equiv 0$}.
The solutions $\zeta_\pm$ are given by
\begin{equation}
\zeta_\pm(z_1,\ldots,z_\kappa) = \half\left[1 + r - \beta(z_1,\ldots,z_\kappa)
     \pm\sqrt{(1 + r - \beta(z_1,\ldots,z_\kappa))^2 - 4\sigma_{K-\kappa}}\right] \;.
\label{ZetaBeta}
\end{equation}
It follows that
\mbox{$\zeta_+() = 1$},
\mbox{$\zeta_-() = r$}.

The explicit expression for the two-level problem
(\mbox{$K = 2$}) is
\begin{align}
\begin{aligned}
G_0(z_1) &= P_0{\cdot}\frac{1 - z_1}{1 - z_1\zeta_+(z_1)} \\
&= P_0{\cdot}\frac{1 - z_1}{1 -\dfrac{z_1}{2}\left[1 + r - r_2z_1 + \sqrt{(1+r-r_2z_1)^2 - 4r_1}\right]} \;,
\end{aligned}
\label{G0K2}
\end{align}
which is in agreement with \citet{NP:Cohen56},
noting that
\mbox{$P_0 = 1-r$}
if $G_0(.)$ is to represent the wait-conditional PGF.
The explicit expression for the three-level problem
(\mbox{$K = 3$}) is
\begin{align}
\begin{aligned}
G_0(z_1,z_2) &= P_0{\cdot}\frac{1 - z_1}{1 - z_1\zeta_+(z_1)}{\cdot}
     \frac{1 - z_2\zeta_+(z_1)}{1 - z_2\zeta_+(z_1,z_2)} \\
&= P_0{\cdot}\frac{1-z_1}{1 - \dfrac{z_1}{2}\left[1 + r - r_3z_1 + \sqrt{(1+r-r_3z_1)^2 - 4(r_1+r_2)}\right]} \\
&\quad {}\times \frac{1 - \dfrac{z_2}{2}\left[1 + r - r_3z_1 + \sqrt{(1+r-r_3z_1)^2 - 4(r_1+r_2)}\right]}
      {1 - \dfrac{z_1}{2}\left[1 + r - r_3z_1 - r_2z_2 + \sqrt{(1+r-r_3z_1 - r_2z_2)^2 - 4r_1}\right]} \;.
\end{aligned}
\label{G0K3}
\end{align}
The result for the full wait-conditional PGF is given by
\begin{equation}
G_\ell(z_1,\ldots,z_{K-1}) = G_0(z_1,\ldots,z_{K-1}){\cdot}\zeta_-^\ell(z_1,\ldots,z_{K-1}) \;,
\label{FullPGF}
\end{equation}
for
\mbox{$\ell = 0,1,\ldots$}.

\subsection{Marginal Probabilities}
We first consider the low-priority marginal PMF for the two-level
(\mbox{$K = 2$})
problem.
The low-priority marginal PGF is a function of a single continuous complex variable $z$,
and is obtained by summing (\ref{FullPGF}) over the discrete high-priority queue lengths, to yield
\begin{equation}
G_{\text{lo}}(z) = \sum_{\ell=0}^\infty G_\ell(z) = \frac{(1-r)(1-z)}
     {[1 - z\zeta_+(z)][1 - \zeta_-(z)]} \;,
\label{GLo1}
\end{equation}
with
\begin{equation}
\zeta_\pm(z) = \half\left[1 + r - r_2z \pm \sqrt{(1 + r - r_2z)^2 - 4r_1}\right] \;,
\label{ZetaK2}
\end{equation}
which agrees with Cohen's result \citep{NP:Cohen56}.
By using the quadratic identities
\begin{equation}
\zeta_+(z) + \zeta_-(z) = 1 + r - r_2z \;, \quad
\zeta_+(z)\cdot\zeta_-(z) = r_1 \;,
\label{ZetaIds}
\end{equation}
which also imply that
\begin{equation}
[1 - \zeta_+(z)]{\cdot}[1 - \zeta_-(z)] = r_2(z-1) \;,
\end{equation}
one is able to express $G_{\text{lo}}(z)$
in various equivalent and more convenient forms:
\begin{equation}
G_{\text{lo}}(z) =  \frac{1-r}{r_2}{\cdot}\frac{r - \zeta_-(z)}{1 - rz}
     = \frac{1-r}{\zeta_+(z) - r} \;.
\label{GLo2}
\end{equation}
Since
\mbox{$\zeta_+(1/r) = r$},
\mbox{$\zeta_-(1/r) = r_1/r$},
the first form shows that $G_{\text{lo}}(z)$ has a pole at
\mbox{$z = 1/r > 1$}.
We also have
\mbox{$\zeta_+(1) = 1$},
\mbox{$\zeta_-(1) = r_1$},
but it should be clear that the apparent pole at
\mbox{$z = 1$}
implied by (\ref{GLo1}) is spurious due to a  cancellation with the numerator.

The $p$-th marginal for the general $K$-level problem, with
\mbox{$p = 1,2,\ldots,K-1$},
has its PGF  given by
\begin{equation}
G_{\text{mrg}}^{(p)}(z_p) \equiv \sum_{\ell=0}^\infty G_\ell(1,\ldots,z_p,\ldots,1)  =
     \frac{G_0(1,\ldots,z_p,\ldots,1)}{1 - \zeta_-(1,\ldots,z_p,\ldots,1)} \;,
\label{PGFzp}
\end{equation}
and pertains to the $p$-th {\em lowest} priority level or, equivalently, the
$(K{+}1{-}p)$-th {\em highest} priority level.
Each $G_{\text{mrg}}^{(p)}(z_p)$ has the same functional form as $G_{\text{lo}}(z)$
for the two-level problem.
Let us relabel (\ref{ZetaK2}) as
\begin{equation}
\zeta_\pm(z; r_{\text{hi}}, r_{\text{lo}}) = \half\left[1 + r_{\text{sum}} -
     r_{\text{lo}}z \pm \sqrt{(1 + r_{\text{sum}} - r_{\text{lo}}z)^2 - 4r_{\text{hi}}}\right] \;,
\end{equation}
by setting
\mbox{$r_{\text{hi}} = r_1$},
\mbox{$r_{\text{lo}} = r_2$},
\mbox{$r_{\text{sum}} \equiv r_{\text{hi}}+r_{\text{lo}}$},
so that we also have
\begin{equation}
G_{\text{lo}}(z; r_{\text{hi}}, r_{\text{lo}}) =
     \frac{1-r_{\text{sum}}}{\zeta_+(z; r_{\text{hi}}, r_{\text{lo}}) - r_{\text{sum}}} \;.
\end{equation}
Then
\mbox{$G_{\text{mrg}}^{(p)}(z_p) = G_{\text{lo}}(z; r_{\text{hi}}, r_{\text{lo}})$}
on making the identifications
\begin{equation}
r_{\text{lo}} = r_{K+1-p}\;, \quad r_{\text{hi}} = \sum_{\kappa=1}^{K-p} r_\kappa \;.
\end{equation}
The derivation of this result is relegated to the second appendix.
Its interpretation is as follows:
All priority levels above the $p$-th lowest level are aggregated into a single high level
whose level traffic intensity is the sum of the aggregated values.
The low level is identified with the $p$-th lowest level and assigned its level traffic intensity,
namely $r_{K+1-p}$. Finally, all priority levels below the $p$-th are
discarded, so the the total traffic intensity for the equivalent two-level problem becomes
\mbox{$r_{\text{sum}} = r_{\text{lo}}+r_{\text{hi}}$}.

\subsection{Exclusively-Low Probabilities}
To consider the case where there are assumed to be only the lowest-priority
clients in the queue, we must set
\mbox{$z_2 = z_3 = \cdots = z_{K-1} = 0$},
to obtain the PGF for the exclusively-low case
\begin{equation}
G_{\text{xlo}}(z_1) = G_0(z_1,0,\ldots,0) \;,
\end{equation}
which is given by (\ref{G0}) and (\ref{FullPGF}) as
\begin{equation}
G_{\text{xlo}}(z_1) = P_0{\cdot}\frac{1-z_1}{1 - z_1\zeta_+(z_1,0,\ldots,0)} \;.
\end{equation}
On setting
\mbox{$r_{\text{lo}} = r_K$},
\mbox{$r_{\text{hi}} = \sigma_{K-1} = r - r_{\text{lo}}$},
we can write
\begin{align}
\begin{aligned}
\beta(z_1,\mathbf{0})   &= r_{\text{lo}}z_1 \;, \\
\zeta_\pm(z_1,\mathbf{0}) &= \half\left[1 + r  - r_{\text{lo}}z_1 \pm
     \sqrt{(1+r-r_{\text{lo}}z_1)^2 - 4r_{\text{hi}}}\right] \;.
\end{aligned}
\end{align}
Thus, we see that
\mbox{$\zeta_\pm(z_1,\mathbf{0})$}
is identical with its counterpart for the two-level ($K = 2$) problem,
so that we can write
\mbox{$\zeta_\pm(z_1,\mathbf{0}) = \zeta_\pm(z_1)$}.
Using the identity
\begin{equation}
\frac{1-z_1}{1 - z_1\zeta_+(z_1)}  = \frac{1-\zeta_-(z_1)}{\zeta_+(z_1)-r} \;,
\end{equation}
we arrive at the result
\begin{equation}
G_{\text{xlo}}(z_1) = P_0{\cdot}\frac{1 - \zeta_-(z_1)}{\zeta_+(z_1) - r}
     = P_0{\cdot}\left[1 + r_{\text{lo}}{\cdot}\frac{z_1}{\zeta_+(z_1) - r}\right] \;.
\end{equation}
This is equivalent to
\begin{equation}
G_{\text{xlo}}(z_1) = 1 - r + r_{\text{lo}}z_1 G_{\text{lo}}(z_1) \;.
\end{equation}
It follows that
\begin{equation}
P_{\text{xlo}}(n) = (1-r)\delta_{n0} + (1 - \delta_{n0})r_{\text{lo}}P_{\text{lo}}(n-1) \;,
\label{Pxlo}
\end{equation}
where
\mbox{$P_{\text{lo}}(n)$}
is the wait-conditional marginal PMF for the lowest-priority level.
One should note that this relationship is independent of the  number of priority levels $K$.

\subsection{Exclusively-High Probabilities}
To consider the case where there are assumed to be only the highest-priority
clients in the queue, we must set
\mbox{$z_1 = z_1 = \cdots = z_{K-1} = 0$},
to obtain the PMF for the exclusively-high case
\begin{equation}
P_{\text{xhi}}(\ell) = G_\ell(0,0,\ldots,0)  = P_0\zeta_-^\ell(\mathbf{0}) \;,
\end{equation}
as given by (\ref{G0}) and (\ref{FullPGF}).
On setting
\mbox{$r_{\text{hi}} = \sigma_1 = r_1$},
\mbox{$r_{\text{lo}} = r - r_{\text{hi}}$},
we can write
\begin{align}
\begin{aligned}
\beta(\mathbf{0})   &= 0 \;, \\
\zeta_\pm(\mathbf{0}) &= \half\left[1 + r  \pm  \sqrt{(1+r)^2 - 4r_{\text{hi}}}\right] \;.
\end{aligned}
\end{align}
Therefore
\begin{equation}
P_{\text{xhi}}(\ell) = (1-r){\cdot}\left[\frac{1 + r - \sqrt{(1+r)^2 - 4r_{\text{hi}}}}{2}\right]^\ell \;,
\label{Pxhi}
\end{equation}
for
\mbox{$\ell = 0,1,\ldots$}.
One should note that this result is independent of the  number of priority levels $K$.

\section{FFT Mixture Method}
\label{FFT}
\begin{figure}
\FIGURE
{\includegraphics[width=\wscl\linewidth, height=\hscl\linewidth]{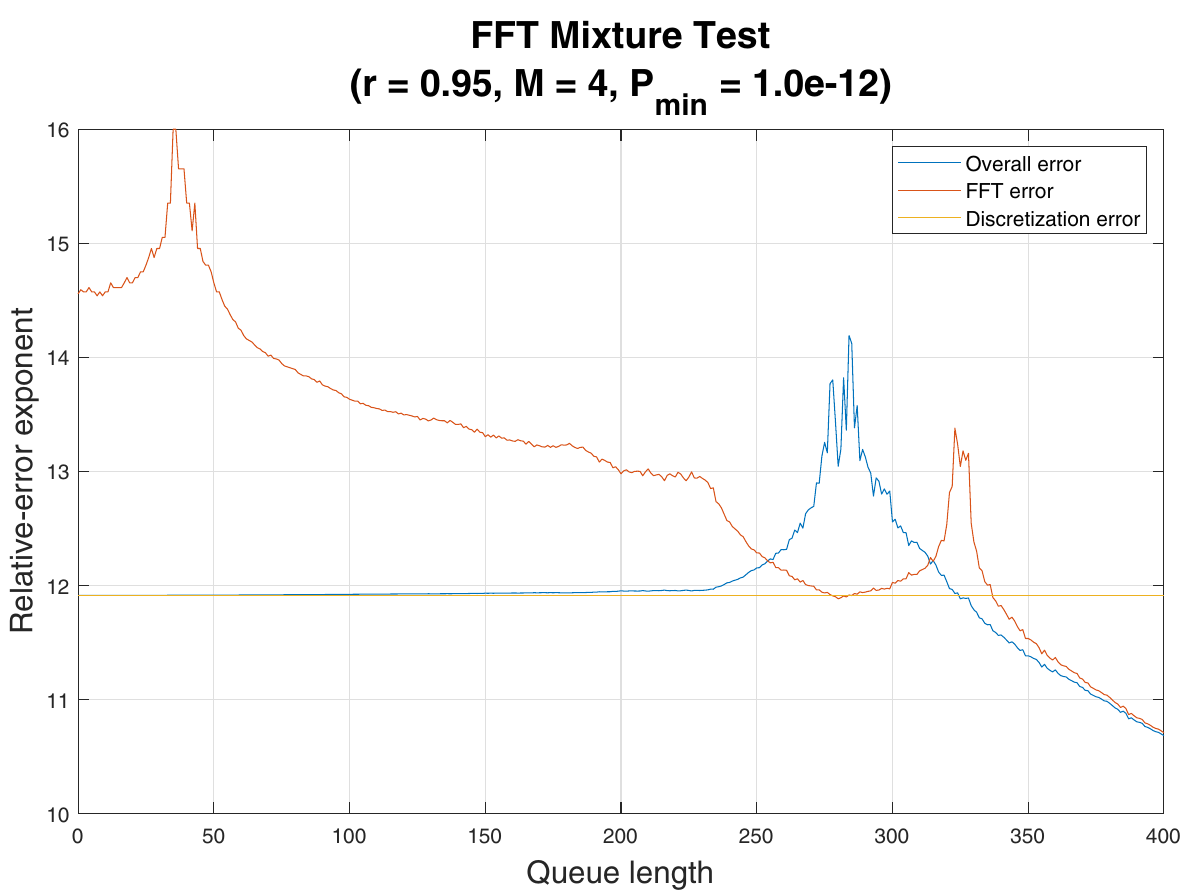}}
{\hphantom{x}\label{FFTMixTestF}}
{Error analysis for the FFT mixture method.}
\end{figure}

Let us begin by considering the numerical evaluation of the marginal PMFs.
Without loss of generality, we may restrict our attention to the low-priority PMF
for the two-level problem.
We have
\begin{equation}
P_{\text{lo}}(n) = \frac{1}{n!}\left.\frac{d^n}{dz^n}g(z)\right|_{z=0} \;, \quad
     g(z) = \frac{1-r}{\zeta_+(z) - r} \;,
\label{gz}
\end{equation}
where
\mbox{$g(z) = G_{\text{lo}}(z)$}
as given by (\ref{GLo2}).
Invoking Cauchy's integral theorem, we can write this as
\begin{equation}
P_{\text{lo}}(n) = \oint_\mathcal{C}\frac{dz}{2\pi i}\, \frac{g(z)}{z^{n+1}} \;,
\end{equation}
for
\mbox{$n = 0,1,\ldots$},
where the integration contour $\mathcal{C}$ encircles the origin anti-clockwise and remains within
the radius of convergence of the Taylor expansion of $g(z)$,
which is given by
\mbox{$\eta_{\text{c}} = 1/r > 1$}.
If we choose $\mathcal{C}$ to be the circle of radius
\mbox{$\eta < \eta_{\text{c}}$}
centred about the origin, and write
\mbox{$z = \eta e^{-i\theta}$},
for
\mbox{$0 \leq \theta < 2\pi$},
then
\begin{equation}
P_{\text{lo}}(n) = \frac{1}{\eta^n}\int_0^{2\pi}\frac{d\theta}{2\pi}\
     e^{in\theta}g\left(\eta e^{-i\theta}\right) \;.
\label{Cauchy}
\end{equation}
Approximation of the integral by an $N$-interval trapezoidal rule on the grid
\mbox{$\theta_k = 2\pi k/N$},
\mbox{$k = 0,1,\ldots,N$},
yields the finite sum
\begin{equation}
P_{\text{lo}}(n) \simeq \frac{1}{\eta^n N}\sum_{k=0}^{N-1} e^{2\pi i nk/N}g\left(\eta e^{-2\pi i k/N}\right) \;,
\label{Trapz}
\end{equation}
which is proportional to the inverse discrete Fourier transform of the sequence
\mbox{$h(k) = g\left(\eta e^{-2\pi i k/N}\right)$},
\mbox{$k = 0,1,\ldots,N-1$}.
When $N$ is chosen to be a power of $2$, it can be implemented as an inverse fast Fourier transform (IFFT),
in which case we write
\mbox{$P_{\text{lo}}(n) \simeq \eta^{-n}{\cdot}{\sf IFFT}[h](n)$}.
Extraction of queue-length probabilities from a PGF by means of a discrete Fourier transform
has been previously considered by \citet{NP:Daigle89}.
We shall adopt an alternative general method due to \citet{NP:Fornberg81}
that is better suited to the multi-variate problem

The RHS of (\ref{Cauchy}) is independent of the value of the contour radius $\eta$.
This is no longer true in (\ref{Trapz}) as a consequence of the approximation.
The dependence on the integration radius in using an FFT to numerically compute derivatives
of analytic functions has been studied by \citet{NP:Bornemann11}.
He found that there exists an optimal radius that minimizes the error, and that this optimal
value is close to the radius of convergence.

The dependence of the approximated solution on the contour radius can be used to one's advantage,
as originally pointed out by \citet{NP:Fornberg81}.
Consider the numerical differentiation of the analytic function $g(z)$.
According to the FFT method, we must evaluate
\begin{equation}
H(n;\eta) \equiv \eta^{-n}{\cdot}{\sf IFFT}_k[h(k;\eta)](n) \;, \quad
     h(k,\eta) \equiv g\left(\eta e^{-2\pi i k/N}\right) \;,
\end{equation}
for some contour radius.
We may equally evaluate the weighted mixture
\begin{equation}
H_M(n) \equiv \sum_{m=1}^M f_m H(n;\eta_m) \;,
\end{equation}
for some collection of contour radii $\eta_m$ and mixture coefficients $f_m$
that sum to unity
\mbox{$\sum_{m = 1}^M f_m = 1$}.
As discussed in \citep{NP:Fornberg81}, the first $M$ aliasing terms can be cancelled
using (in principle) any $M$ distinct contour radii
\mbox{$\eta_m < \eta_{\text{c}}$}
provided ones makes an appropriate choice of the coefficients $f_m$.
To achieve this, the mixture coefficients $f_m$
must solve the matrix equation
\begin{equation}
\begin{bmatrix}
1               &  1              & \cdots & 1             \\
\eta_1^N        & \eta_2^N        & \cdots & \eta_M^N       \\
\eta_1^{2N}     & \eta_2^{2N}     & \cdots & \eta_M^{2N}    \\
\vdots          & \vdots          & \ddots & \vdots         \\
\eta_1^{(M-1)N} & \eta_2^{(M-1)N} & \cdots & \eta_M^{(M-1)N}
\end{bmatrix} \!\!
\begin{bmatrix}
f_1    \vphantom{1}            \\
f_2    \vphantom{\eta_M}       \\
f_3    \vphantom{\eta_M^2}     \\
\vdots \vphantom{\vdots}       \\
f_M    \vphantom{\eta_M^{M-1}}
\end{bmatrix} =
\begin{bmatrix}
1 \\
0 \\
0 \\
\vdots \\
0
\end{bmatrix} \;.
\label{MatEq}
\end{equation}
The solution for the mixture coefficients $f_m$ is given by
\begin{equation}
\frac{1}{f_m} = \prod_{\substack{\ell = 1\\ \ell\neq m}}^M\left(1 - \frac{\eta_m^N}{\eta_\ell^N}\right) \;,
\label{fm}
\end{equation}
for
\mbox{$m = 1,2,\dots,M$}.
%
Therefore, the residual errors after the FFT mixture are given by
\begin{equation}
\varepsilon_M(n) = c_{M\!N}(n)\sum_{m=1}^M f_m\eta_m^{MN}
     = (-1)^{M-1}c_{M\!N}(n) \prod_{m=1}^M\eta_m^N \;,
\label{ResErr}
\end{equation}
for some constants $c_{M\!N}(n)$.
The result for the summation over $m$ follows directly from (\ref{fm}).
One may note that, when applied to (\ref{gz}),
\mbox{$c_{M\!N}(n) = P_{\text{lo}}(n+M\!N)$}.

\subsection{Marginal Distribution}
Without loss of generality, we consider the marginal distribution for the low-priority level in the two-level
(\mbox{$K=2$}) problem.
For any fixed $r$, worst case performance of the FFT method is observed to occur as
\mbox{$r_{\text{hi}}\to 0^+$}.
Thus, we shall set
\mbox{$r_{\text{lo}} = r$},
\mbox{$r_{\text{hi}} = 0$}.
In this case, the wait-conditional PGF is trivially given by
\begin{equation}
g(z) = \frac{1- r}{1 - rz} = \sum_{n=0}^\infty p_n z^n \;,
\end{equation}
with
\mbox{$p_n = (1-r)r^n$},
for
\mbox{$k = 0,1,\ldots$}.
To approximate the $p_n$ via an $N$-point FFT, we apply the trapezoidal rule on the grid
\mbox{$z_k = \eta e^{-2\pi ik/N}$}
for
\mbox{$k = 0,1,\ldots,N-1$},
to obtain
\begin{equation}
p_n = \frac{1}{N}\sum_{k=0}^{N-1} g(z_k)e^{2\pi ikn/N} \;.
\end{equation}
This leads us to consider the ratio
\begin{equation}
\frac{p_n}{(1-r)r^n} =\frac{1}{N(r\eta)^n}\sum_{k=0}^{N-1}
     \frac{e^{2\pi ikn/N}}{1 - r\eta e^{-2\pi ik/N}}
\label{prat}
\end{equation}
that will approach unity as
\mbox{$N\to\infty$}.
With
\mbox{$\xi \equiv r\eta < 1$},
we introduce the ratio function
\begin{equation}
R_N(\xi,n) \equiv \frac{1}{N\xi^n}\sum_{k=0}^{N-1}\frac{e^{2\pi ikn/N}}
     {1 - \xi e^{-2\pi ik/N}} \;.
\end{equation}
Given (\ref{prat}), we aim to achieve
\mbox{$R_N(\xi,n) \simeq 1$}
for all
\mbox{$n = 0,1,\ldots,N-1$}.
We may also observe that
\begin{equation}
R_N(\xi,n) = \frac{1}{\xi^n}I_N(\xi, n)
     = \frac{1}{\xi^n}\ifft_k\left[(1 - \xi e^{-2\pi ik/N})^{-1}\right](n) \;.
\label{RIN}
\end{equation}
The ratio function can be evaluated by explicitly summing the series to give
\begin{equation}
R_N(\xi,n) = \frac{1}{1 - \xi^N} \;,
\end{equation}
independent of $n$.
Thus,
\mbox{$\delta R_N(\xi) \equiv R_N(\xi,n) - 1 \simeq \xi^N$}
for
\mbox{$\xi^N \ll 1$}.
If the desired discretization relative error is set to $10^{-\alpha}$,
then $\xi$ is determined by
\mbox{$\xi^N = 10^{-\alpha}$}.

We shall be directly computing the quantity $I_N(\xi, n)$ in (\ref{RIN})
via a dedicated FFT routine.
Since
\mbox{$R_N(\xi,n) \simeq 1$},
we have that
\mbox{$I_N(\xi, n) \sim \xi^n$}.
This cannot be too small for relevant values of $n$ due to machine arithmetic limits.
In double-precision arithmetic, we must ensure that
\mbox{$I_N(\xi, n) \gg 10^{-16}$}
for the largest
\mbox{$n = 0,1,\ldots,N_{\text{max}} < N$}
of interest, $N$ being the FFT size,
which must be greater than $N_{\text{max}}$ in order to prevent aliasing.
If we take the absolute FFT error ({\it i.e.}\ round-off error) to be
\mbox{$\varepsilon_{\text{fft}} \simeq 10^{-15}$},
then the relative FFT error is
\mbox{$10^{-15}/\xi^{N_{\text{max}}} \sim 10^{-15 + \alpha N_{\text{max}}/N}$}.
Best numerical performance is obtained by equating the two competing and opposing
sources of error, {\it i.e.}\ the relative discretization error with the relative FFT error.
This yields the relationship
\begin{equation}
\alpha = \frac{15}{1 + N_{\text{max}}/N}
\end{equation}
for the overall relative error exponent.
It implies that, for twelve decimal places of accuracy
(\mbox{$\alpha = 12$}),
we require that
\mbox{$N \simeq 4N_{\text{max}}$}.

Let us now consider the FFT mixture scheme.
In order to specify the $M$ contour radii,
we introduce a spread parameter $s$, whose value we typically take to be
\mbox{$s \simeq 0.05$},
and spread factors
\mbox{$\varsigma_m \equiv 1 - s{\cdot}(m-1)/(M-1)$},
for
\mbox{$m = 1,2,\ldots,M$}.
Then, we set
\mbox{$\xi_m = \varsigma_m\xi$},
which gives all of the radii in terms of the largest one
\mbox{$\xi_1 = \xi$},
as yet undetermined.
The smallest is given by
\mbox{$\xi_M = \varsigma_M \xi = (1-s)\xi$}.
It is useful to observe that, in the present context, (\ref{ResErr}) reads
\begin{equation}
\varepsilon_M(n)/[(1-r)r^n] = (-1)^{M-1}\prod_{m=1}^M\xi_m^N \;,
\end{equation}
the RHS being independent of $n$.

We revisit the foregoing error analysis adapted to the mixture method,
and begin by estimating the FFT error.
If the mixture ratio function as given by the trapezoidal rule were the
actual object of evaluation by means of an FFT, then the exact result would be
\begin{equation}
R^{\text{mix}}_{N_{\text{fft}}}(\boldsymbol{\xi},n) =
     \sum_{m=1}^M\frac{f_m(\boldsymbol{\xi})}{1 - \xi_m^{N_{\text{fft}}}} \;,
\end{equation}
while the approximate result rendered by the FFT is
\begin{align}
\begin{aligned}
R^{\text{fft}}_{N_{\text{fft}}}(\boldsymbol{\xi},n) &=
     \sum_{m=1}^M\frac{f_m(\boldsymbol{\xi})}{1 - \xi_m^{N_{\text{fft}}}}{\cdot}
     \left(1 \pm \frac{\varepsilon_{\text{fft}}}{\xi_m^n}\right) \\
&= R^{\text{mix}}_{N_{\text{fft}}}(\boldsymbol{\xi},n) \pm \varepsilon_{\text{fft}}
     \sum_{m=1}^M\frac{|f_m(\boldsymbol{\xi})|}{1 - \xi_m^{N_{\text{fft}}}}{\cdot}\frac{1}{\xi_m^n} \;,
\end{aligned}
\end{align}
where $\varepsilon_{\text{fft}}$ denotes the absolute error in the FFT computation due to
machine arithmetic precision.
Thus, the error in the ratio function due to FFT precision is
\begin{align}
\begin{aligned}
\Delta R_{N_{\text{fft}}}^{\text{fft}}(\boldsymbol{\xi}, n) &\equiv \left|
     R^{\text{fft}}_{N_{\text{fft}}}(\boldsymbol{\xi},n) -
     R^{\text{mix}}_{N_{\text{fft}}}(\boldsymbol{\xi},n)\right| \\
&= \varepsilon_{\text{fft}}
     \sum_{m=1}^M\frac{f_m(\boldsymbol{\xi})}{1 - \xi_m^{N_{\text{fft}}}}{\cdot}\frac{1}{\xi_m^n} \\
&< \varepsilon_{\text{fft}}
     \sum_{m=1}^M\frac{|f_m(\boldsymbol{\xi})|}{1 - \xi_m^{N_{\text{fft}}}}{\cdot}\frac{1}{\xi_m^n} \;,
\end{aligned}
\end{align}
where the final step represents a worst-case bound.
Given that
\mbox{$\xi_m^N \ll 1$},
we have
\begin{equation}
\sum_{m=1}^M\frac{|f_m(\boldsymbol{\xi})|}{1 - \xi_m^N}{\cdot}\frac{1}{\xi_m^n} \simeq
     \sum_{m=1}^M\frac{|f_m(\boldsymbol{\xi})|}{\xi_m^n} \simeq
     \frac{|f_M(\boldsymbol{\xi})|}{\xi_M^n} \;.
\end{equation}
So, we see that the FFT error is dominated by the smallest contour radius $\xi_{M}$.
It follows that it is estimated by
\begin{equation}
\Delta R_{N_{\text{fft}}}^{\text{fft}}(\boldsymbol{\xi}, n) \sim
     \varepsilon_{\text{fft}}{\cdot}|f_M(\boldsymbol{\xi})|/\xi_M^n \sim
     \varepsilon_{\text{fft}}/\xi_M^n \;,
\label{DRfft}
\end{equation}
since
\mbox{$|f_M(\boldsymbol{\xi})| \simeq 1$}
for a sufficiently large spread $s$, and we take
\mbox{$\varepsilon_{\text{fft}} = 10^{-15}$}.

To estimate the discretization error, we consider the mixture ratio function
\begin{align}
\begin{aligned}
R^{\text{mix}}_{N_{\text{fft}}}(\boldsymbol{\xi},n) &=
     \sum_{m=1}^M\frac{f_m(\boldsymbol{\xi})}{1 - \xi_m^{N_{\text{fft}}}} \\
&= \sum_{m=1}^M f_m(\boldsymbol{\xi}) + \sum_{m=1}^M \xi_m^{N_{\text{fft}}} f_m(\boldsymbol{\xi}) + \cdots \\
&= 1 + (-1)^{M-1}(\xi_1\xi_2\cdots\xi_M)^{N_{\text{fft}}} + \cdots \;.
\end{aligned}
\end{align}
Thus, the discretization error due to application of the trapezoidal rule is
\begin{equation}
\Delta R^{\text{tpz}}_{N_{\text{fft}}}(\boldsymbol{\xi},n) \equiv
     \left|R^{\text{mix}}_{N_{\text{fft}}}(\boldsymbol{\xi},n)- 1\right| \simeq
     (g\xi)^{MN_{\text{fft}}} \;,
\label{DRtpz}
\end{equation}
where we have introduced the geometric mean
\begin{equation}
g \equiv \biggl(\prod_{m=1}^M\varsigma_m\biggr)^{1/M} \;.
\end{equation}
We proceed to equate the trapezoidal-rule error of (\ref{DRtpz})
with the FFT error of (\ref{DRfft}) at the largest desired queue length
$n = N_{\text{max}}$,
and set them to a common error level
$10^{-\alpha}$, {\it i.e.}\
\begin{equation}
(g\xi)^{MN_{\text{fft}}} = 10^{-15}/\xi_M^{N_{\text{max}}} = 10^{-\alpha} \;.
\end{equation}
Then, after taking
\mbox{$\alpha = 12$},
\mbox{$M = 4$},
eliminating $\xi$ from the equations yields the relationship
\mbox{$N_{\text{fft}} = N_{\text{max}}/(1 - N_{\text{max}}\chi)$}
with
\mbox{$\chi \equiv \log_{10}((1-s)/g)$}.
For
\mbox{$s = 0.05$,}
\mbox{$M = 4$},
we obtain
\mbox{$\chi = 0.0037$}.
So, it suffices to choose
\mbox{$N_{\text{fft}} \simeq N_{\text{max}}$}.
Consequently, the largest contour radius is determined as
\begin{equation}
\xi = g^{-1}{\cdot}10^{-12/(MN_{\text{fft}})} \;,
\end{equation}
with $N_{\text{fft}}$ chosen to be the smallest power of two exceeding $N_{\text{max}}$.

As a sanity check, we apply the FFT-mixture method as described here to the computation
of the low-priority marginal PMF for the two-level problem of the case
\mbox{$r_{\text{hi}} = 0$},
and plot in Figure~\ref{FFTMixTestF} various errors for the ratio function
as a function of queue length.
The maximum queue length considered was that sufficient to attain a tail probability
\mbox{$P_{\text{min}} = 10^{-12}$}.
The target error level was set to $10^{-12}$
({\it i.e.}\ $\alpha = 12$).
The exact (pre-discretization) ratio function is everywhere unity in the present case.
The overall error (blue curve) represents the difference in the computed ratio function from unity.
The FFT error (red curve) represents the difference in the computed ratio function from the exact value
of the discretized ratio function.
The discretization error (orange curve) represents the difference in the exact discretized ratio function
from unity.
The results are consistent with expectations.

\subsection{Joint Distribution}
The expression (\ref{G0}) for $G_0(\mathbf{z})$ contains spurious algebraic singularities that cancel
out between denominator and numerator, as can be seen by observing that
\mbox{$\zeta_+(\mathbf{1}) = 1$}.
This has implications for numerical evaluation, and limits the choice of integration contour radii
in the application of Cauchy's theorem to
\mbox{$|z_\kappa| < 1$}
whereas, in principle, one could extend this to at least
\mbox{$|z_\kappa| < 1/r$}.
Fortunately, the algebraic singularities can be eliminated by re-structuring the integrand
into an equivalent form according to the following argument:
Let us write (\ref{G0}) as
\begin{equation}
G_0(\mathbf{z}) = P_0{\cdot}\prod_{\kappa=1}^{K-1}
     \mathscr{P}^+_\kappa(\mathbf{z})/\mathscr{Q}^+_\kappa(\mathbf{z}) \;,
\label{G0PQ}
\end{equation}
where
\begin{equation}
\mathscr{P}^\pm_\kappa(\mathbf{z}) \equiv 1 - z_\kappa\zeta^\pm_{\kappa-1}(\mathbf{z}) \;, \quad
\mathscr{Q}^\pm_\kappa(\mathbf{z}) \equiv 1 - z_\kappa\zeta^\pm_\kappa(\mathbf{z}) \;,
\label{PQ}
\end{equation}
with
\mbox{$\zeta^\pm_\kappa(\mathbf{z}) \equiv \zeta_\pm(z_1,\ldots,z_\kappa)$},
for
\mbox{$\kappa = 0,1,\ldots,K-1$}.
In particular,
\mbox{$\zeta^+_0() = 1$},
\mbox{$\zeta^-_0() = r$}.
Then, as derived in Appendix \ref{AppRatio}, we have the result
\begin{equation}
\mathscr{P}^+_\kappa(\mathbf{z})/\mathscr{Q}^+_\kappa(\mathbf{z}) =
     \mathscr{Q}^-_\kappa(\mathbf{z})/\mathscr{P}^-_\kappa(\mathbf{z}) \;,
\end{equation}
which yields the equivalent form for the PGF
\begin{equation}
G_\ell(z_1,\ldots,z_{K-1}) = P_0{\cdot}\prod_{\kappa=1}^{K-1}\frac{1 - z_\kappa\zeta_-(z_1,\ldots,z_{\kappa})}
     {1 - z_\kappa\zeta_-(z_1,\ldots,z_{\kappa-1})}{\cdot}\zeta_-^\ell(z_1,\ldots,z_{K-1}) \;.
\label{GMinus}
\end{equation}
This result can also be established directly via an argument based on the fact that
the marginal queue-length distribution for the aggregation of the highest $p$
priority levels is a geometric distribution, for all
\mbox{$p = 1,2,\ldots,K$}.
The details are presented in Appendix~\ref{AltMulti}.
It is also interesting to note that the change of integration variables
\mbox{$\mathbf{z} \mapsto \mathbf{w}$}
such that
\mbox{$w_\kappa(z_1,\ldots,z_\kappa) = \zeta_-(z_1,\ldots,z_\kappa)$},
for
\mbox{$\kappa = 1,2\ldots,K-1$},
leads to the completely meromorphic form
\begin{equation}
G_\ell(z_1,\ldots,z_{K-1}){\cdot}\prod_{\kappa=1}^{K-1}dz_\kappa =
     P_0 w^\ell_{K-1}{\cdot}\prod_{\kappa=1}^{K-1}
     \frac{1 - z_\kappa w_\kappa}{1 - z_\kappa w_{\kappa-1}}{\cdot}\prod_{\kappa=1}^{K-1}
     \frac{\sigma_{K-\kappa} - w^2_\kappa}{r_{\kappa+1}w^2_\kappa}dw_\kappa \;,
\label{GMerom}
\end{equation}
where
\begin{align}
\begin{aligned}
z_\kappa &= \frac{w_\kappa(w^2_{\kappa-1} + \sigma_{K+1-\kappa}) - \sigma_{K-\kappa}w_{\kappa-1}
     - w_{\kappa-1}w^2_\kappa}{r_{K+1-\kappa}w_{\kappa-1}w_\kappa} \\
&= \frac{1}{r_{K+1-\kappa}}\left(w_{\kappa-1} + \frac{\sigma_{K+1-\kappa}}{w_{\kappa-1}}
     - w_\kappa - \frac{\sigma_{K-\kappa}}{w_\kappa}\right) \;,
\end{aligned}
\end{align}
and we adopt the convention that
\mbox{$z_0 \equiv r = \sum_{\kappa=1}^K r_\kappa$}.
The second product is the Jacobian of the variable transformation.
We have not yet found a practical use for this representation.

Application of Cauchy's integral theorem in each of $K{-}1$ dimensions yields the
expression for the joint PMF
\begin{equation}
P(\ell,n_{K-1},\ldots,n_1) = \prod_{\kappa=1}^{K-1}\biggl\{\oint_{\mathcal{C}_\kappa}
     \frac{dz_\kappa}{2\pi iz_\kappa^{n_\kappa+1}}\biggr\}\,
     {\cdot}G_\ell(z_1,\ldots,z_{K-1}) \;,
\end{equation}
where the anti-clockwise closed contours around the origin $\mathcal{C}_\kappa$ can be taken
to be circles of radius
\mbox{$\eta_\kappa < 1/r$}.
For simplicity, we adopt a common radius
\mbox{$\eta_\kappa = \eta$}
in each dimension.

We introduce the multi-dimensional inverse FFT
\begin{equation}
\mathscr{F}_\ell(\mathbf{n};\eta) = \eta^{-N(K-1)}{\cdot}
      {\sf IFFT}_\mathbf{k}[G_\ell(\eta e^{-2\pi i\mathbf{k}/N})](\mathbf{n}) \;,
\end{equation}
where, for simplicity, the FFT size $N$ is taken to be the same for every dimension, and
\mbox{$\mathbf{n} = [n_1,n_2,\ldots,n_{K-1}]$}.
According to the mixture method,
\begin{equation}
P(\ell,n_{K-1},\ldots,n_1) \simeq \sum_{m=1}^M f_m \mathscr{F}_\ell(\mathbf{n};\eta_m) \;,
\end{equation}
with the contour radii $\eta_m$ and the coefficients $f_m$ chosen as for the
marginal distributions discussed in the preceding section.

Using the mixture, we execute $M$ runs of an FFT of size $N$ in each of the $K{-}1$ dimensions.
Without the mixture technique,
comparable accuracy would require a single FFT execution of size $MN$.
Therefore the no-mixture/mixture timing ratio for $K$ priority levels is
\begin{equation}
\frac{(MN)^{K-1}\log_2(MN)^{K-1}}{M{\cdot}N^{K-1}\log_2(N^{K-1})} =
     M^{K-2}{\cdot}\left[1 + \frac{\log_2 M}{\log_2 N}\right] \simeq M^{K-2} \;.
\end{equation}
We see that there is no benefit in implementing the mixture method for the
two-level
(\mbox{$K=2$})
problem, but it is increasingly advantageous as the number of priority levels rises.
Testing of numerous cases confirms this performance differential.

\section{Numerical Tests}
\label{Tests}
\subsection{Aggregation Test}
As discussed for the two-level problem in \citep{NP:Zuk23},
the aggregated queue-length distribution describes the total number of entities
in the queue, regardless of priority level.
This is equivalent to the queue-length distribution of the basic M/M/$c$ queueing model
with traffic intensity
\mbox{$r = \sum_{\kappa=1}^K r_\kappa$},
which is known to be a simple geometric distribution.
Hence, the exact aggregate PMF is given by
\begin{equation}
P_{\text{agg}}^{(\text{ex})}(k) = (1-r)r^k \;,
\label{aggtest}
\end{equation}
for
\mbox{$k = 0,1,2,\ldots$}.
Consequently, for
\mbox{$k > 0$},
\begin{equation}
[\Delta\ln P_{\text{agg}}](k) \equiv \ln(P_{\text{agg}}(k)) - \ln(P_{\text{agg}}(k-1))  = \ln r \;,
\end{equation}
independent of $k$.

One diagnostic test of the computational methodology is to check how well
the aggregate PMF constructed from the computed joint PMF reproduces the exact result.
This test is more convenient than similarly testing against marginals as only
finite summations are required.
The aggregate PMF is obtained from the joint PMF as
\begin{equation}
P_{\text{agg}}(k) = \sum_{\substack{\mathbf{n}\in\mathbb{N}_0^K \\ n_1+\cdots+n_K = k}} P(\mathbf{n}) \;,
\end{equation}
for
\mbox{$k = 0,1,2,\ldots$}.
We then consider the measure of performance (MOP)
\begin{equation}
\Xi_{\text{agg}} \equiv -\max_{k \geq 1}\left\{\log_{10}\left(|[\Delta\ln P_{\text{agg}}](k) -
    \Delta\ln P_{\text{agg}}^{(\text{ex})}|\right)\right\} \;,
\end{equation}
where the maximum is taken over all values
\mbox{$0 \leq k \leq n_{\text{lim}}$}
such that
\mbox{$P_{\text{agg}}^{(\text{ex})}(k) > P_{\text{min}} > 0$},
for some threshold level $P_{\text{min}}$ since one cannot expect the
numerical methods to maintain performance down to
arbitrarily small tail probabilities.
Since we are working in double-precision arithmetic, all MOPs of this kind are capped
at a maximum allowed value of $16$.
The interpretation of $\Xi_{\text{agg}}$
(and similarly for all of the subsequent MOPs)
is that it indicates the number of decimal places of numerical agreement
in the worst case.

\subsection{Nearest-Neighbour Test}
A direct consequence of the identity (\ref{Interior}) is that
the joint PMF at any given interior point
\mbox{$\mathbf{n}$}
is a positively weighted sum of the joint PMF values at $K+1$
of its $2K$ nearest neighbours.
We denote by
\mbox{$P_{\text{nn}}(\mathbf{n})$}
the PMF evaluated at the point
\mbox{$\mathbf{n}$}
by means of (\ref{Interior}).
Then, we consider the MOP
\begin{equation}
\Xi_{\text{nn}} \equiv -\max_{\mathbf{n}\in\mathbb{N}^K}\left\{\log_{10}\left(|\ln(P(\mathbf{n})) -
    \ln(P_{\text{nn}}(\mathbf{n}))|\right)\right\} \;,
\end{equation}
where the maximum is taken over all values
\mbox{$\mathbf{n}\in [1,n_{\text{lim}}]^K$}
such that
\mbox{$P(\mathbf{n}) > P_{\text{min}} > 0$}.

\subsection{Xhi-Test}
We denote by
\mbox{$P_{\text{xhi}}^{(\text{ex})}(\ell)$}
the theoretical exact expression of the exclusively-high PMF as given by (\ref{Pxhi}),
which we compare with the result obtained directly from the joint PMF as
\begin{equation}
P_{\text{xhi}}(\ell) = P(\ell,0,\ldots,0) \;.
\end{equation}
If we define, for all
\mbox{$\ell > 0$},
\begin{equation}
[\Delta\ln P_{\text{xhi}}](\ell) \equiv \ln(P_{\text{xhi}}(\ell)) - \ln(P_{\text{xhi}}(\ell-1)) \;,
\end{equation}
then
\mbox{$\Delta\ln P_{\text{xhi}}^{(\text{ex})} =  \ln\bigl([1+ r - \sqrt{(1+r)^2 - 4r_{\text{hi}}}]/2\bigr)$},
independent of $\ell$.
The MOP for the xhi-test is then taken to be
\begin{equation}
    \Xi_{\text{xhi}} \equiv -\max_{\ell \geq 1}\left\{\log_{10}\left(|[\Delta\ln P_{\text{xhi}}](\ell) -
    \Delta\ln P_{\text{xhi}}^{(\text{ex})}|\right)\right\} \;,
\end{equation}
where the maximum is taken over all values
\mbox{$0 \leq \ell \leq n_{\text{lim}}$}
such that
\mbox{$P_{\text{xhi}}^{(\text{ex})}(\ell) > P_{\text{min}} > 0$}.
When applied to the FFT method, the xhi-test simply checks how accurately
the numerical FFT mixture algorithm
reproduces the Cauchy residue theorem for a simple pole at the origin.

\subsection{Xlo-Test}
The exclusively-low PMF, as obtained directly from the computed joint PMF,
is given by
\begin{equation}
P_{\text{xlo}}(n) = P(0,\ldots,0,n) \;.
\end{equation}
It is theoretically related to the lowest-priority marginal PMF
according to (\ref{Pxlo}).
This allows one to compare the computed $P_{\text{xlo}}(n)$ against
a known accurate evaluation of the marginal $P_{\text{lo}}(n)$.
Such a marginal PMF is provided by the quadratic recurrence method
developed in \citep{NP:Zuk23} for the two-level problem,
recalling that the two-level problem is sufficient for the
calculation of all marginals of the multi-level problem.
The relevant MOP is taken to be
\begin{equation}
    \Xi_{\text{xlo}} \equiv -\max_{n > 0}\left\{\log_{10}\left(|\ln(P_{\text{xlo}}(n)) -
    \ln{(}r_{\text{lo}}P_{\text{lo}}(n-1))|\right)\right\} \;,
\end{equation}
where the maximum is taken over all values
\mbox{$0 < n \leq n_{\text{lim}}$}
such that
\mbox{$P_{\text{lo}}(n) > P_{\text{min}} > 0$}.

\subsection{FPI Test}
In the FPI test, we compute the joint PMF using both the FPI and FFT methods,
and check their agreement. The MOP is simply taken to be
\begin{equation}
    \Xi_{\text{fpi}} \equiv -\max_{\mathbf{n}\in\mathbb{N}^K}\left\{\log_{10}
    \left(|\ln(P_{\text{fft}}(\mathbf{n})) -
    \ln(P_{\text{fpi}}(\mathbf{n}))|\right)\right\} \;,
\end{equation}
where the maximum is taken over all values
\mbox{$\mathbf{n}\in [1,n_{\text{lim}}]^K$}
such that
\mbox{$P_{\text{fft}}(\mathbf{n}) > P_{\text{min}} > 0$}.

\begin{figure}
\FIGURE
{\includegraphics[width=\wscl\linewidth, height=\hscl\linewidth]{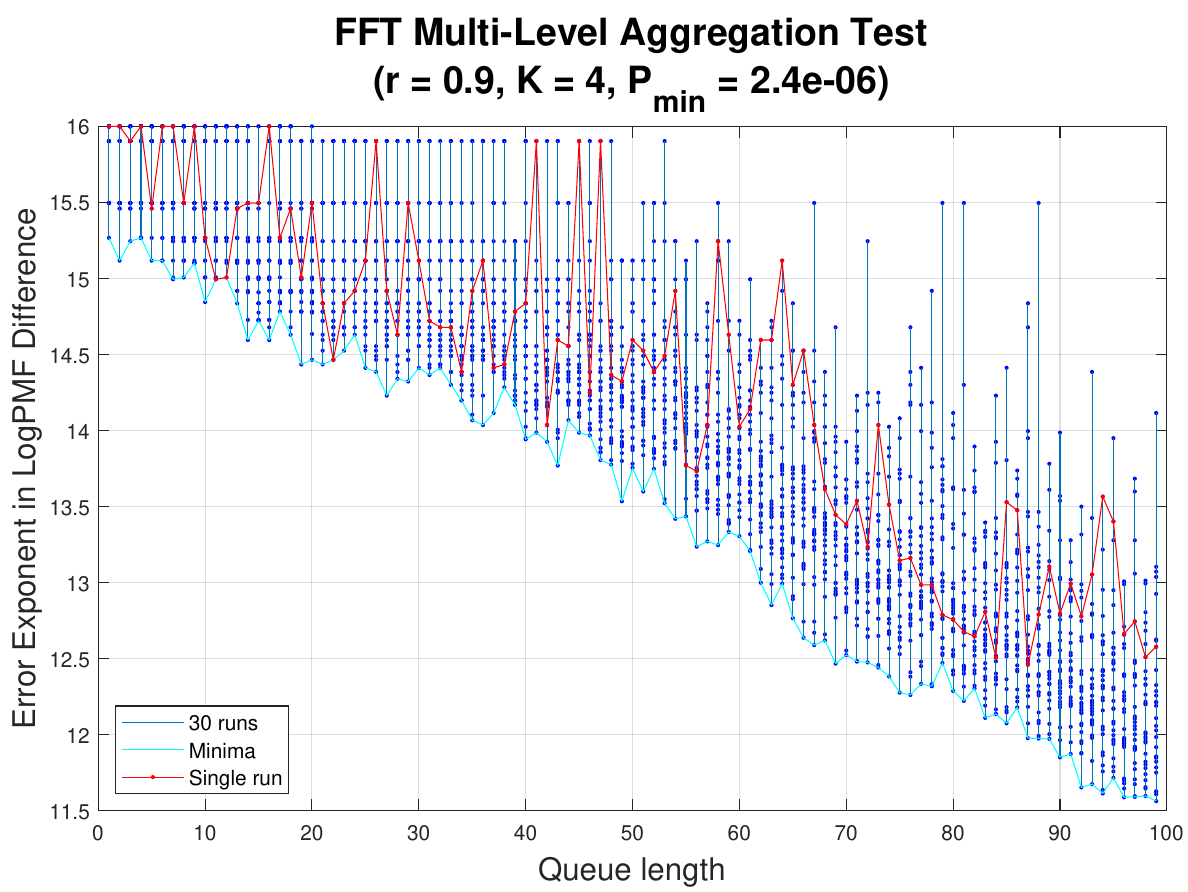}}
{\hphantom{x}\label{TestRepAggFFTJointLevelPMF}}
{Aggregation test for $K = 4$ priority levels with $30$ randomized level traffic intensities corresponding to
     total traffic intensity $r = 0.9$. The number of decimal places of agreement with the exact result is
     plotted on the vertical axis as a function of aggregated queue length. The maximum aggregated queue length of
     $100$ includes data points with PMF above $P_{\text{min}} = 2.4\times 10^{-6}$.}
\end{figure}

\begin{figure}
\FIGURE
{\includegraphics[width=\wscl\linewidth, height=\hscl\linewidth]{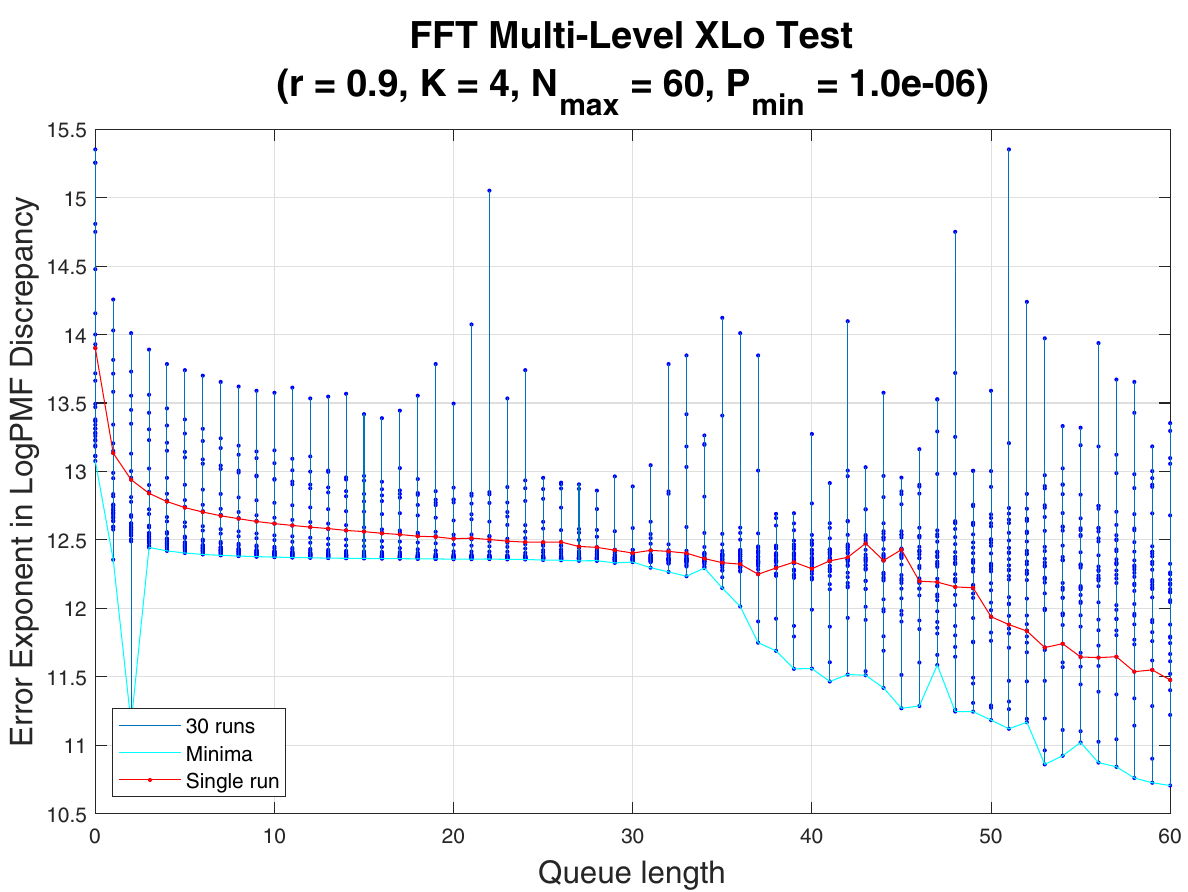}}
{\hphantom{x}\label{TestRepXLoFFTJointLevelPMF}}
{Xlo-test for $K = 4$ priority levels with $30$ randomized level traffic intensities corresponding to
     total traffic intensity $r = 0.9$. The number of decimal places of agreement is
     plotted on the vertical axis as a function of lowest-priority queue length.
     The maximum queue length of
     $60$ includes data points with PMF above $P_{\text{min}} = 1.0\times 10^{-6}$.}
\end{figure}

\begin{figure}
\FIGURE
{\includegraphics[width=\wscl\linewidth, height=\hscl\linewidth]{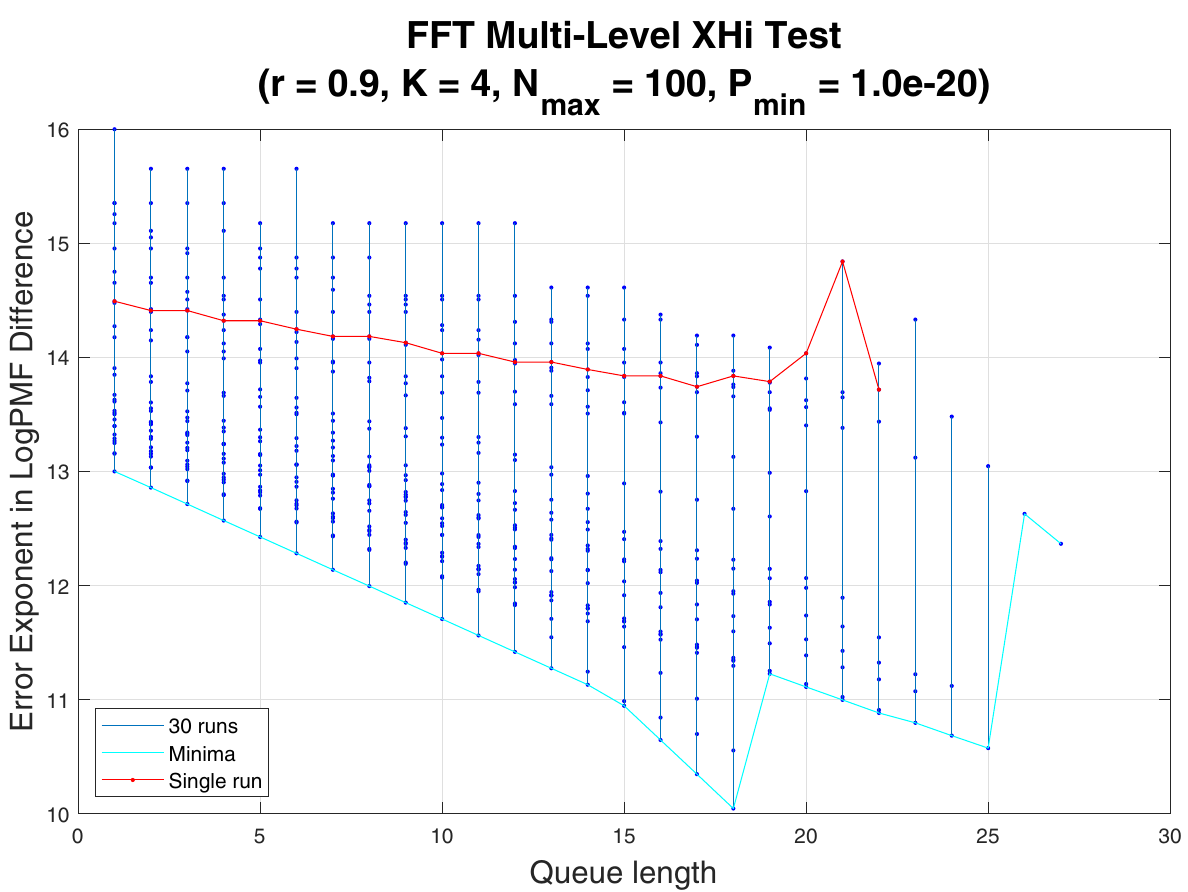}}
{\hphantom{x}\label{TestRepXHiFFTJointLevelPMF}}
{Xhi-test for $K = 4$ priority levels with $30$ randomized level traffic intensities corresponding to
     total traffic intensity $r = 0.9$. The number of decimal places of agreement is
     plotted on the vertical axis as a function of highest-priority queue length.
     The maximum queue length of
     $100$ includes data points with PMF above $P_{\text{min}} = 1.0\times 10^{-20}$.}
\end{figure}

\begin{figure}
\FIGURE
{\includegraphics[width=\wscl\linewidth, height=\hscl\linewidth]{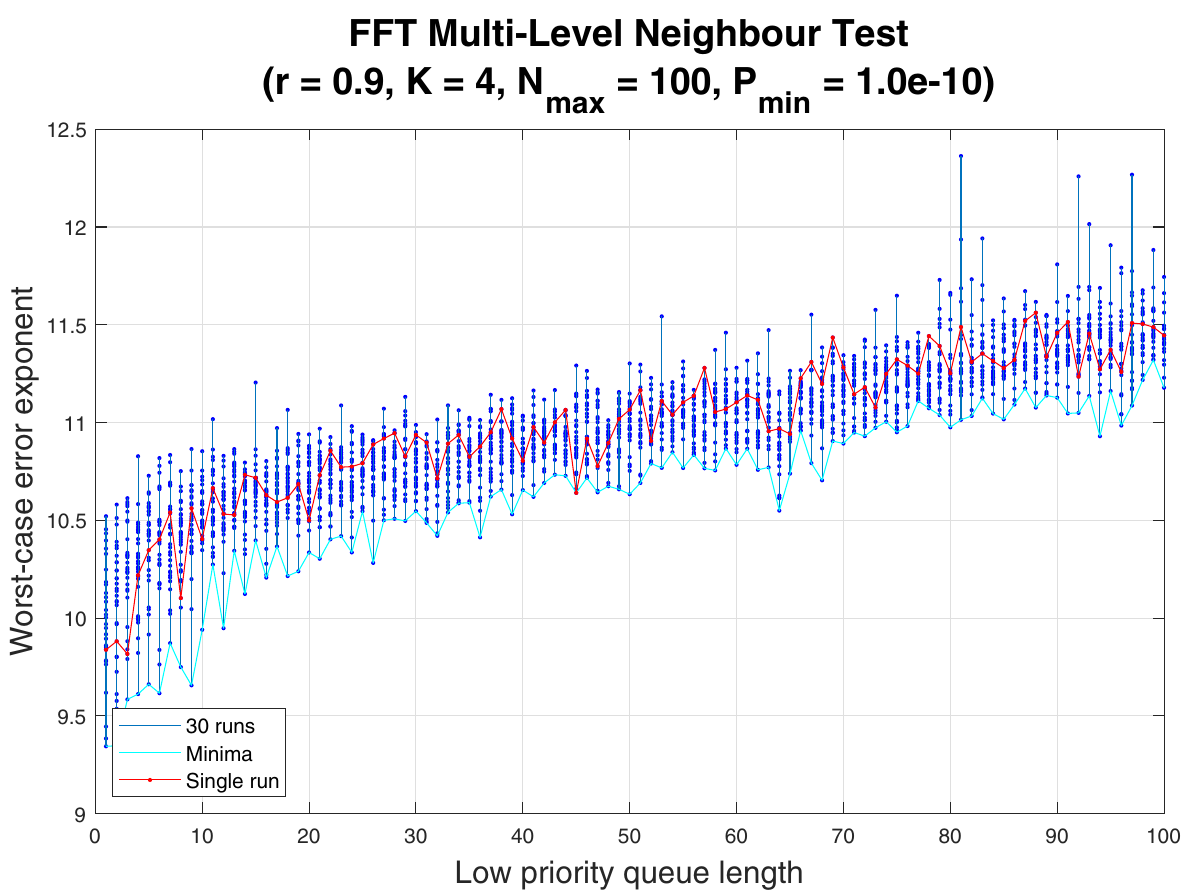}}
{\hphantom{x}\label{TestRepNNFFTJointLevelPMF}}
{Nearest-neighbour test for $K = 4$ priority levels with $30$ randomized level traffic intensities corresponding to
     total traffic intensity $r = 0.9$. The number of decimal places of agreement is
     plotted on the vertical axis as a function of lowest-priority queue length.
     All joint-PMF data points above $P_{\text{min}} = 1.0\times 10^{-10}$ that occur within a
     maximum queue length of $100$ in each priority dimension were considered.}
\end{figure}

\begin{figure}
\FIGURE
{\includegraphics[width=\wscl\linewidth, height=\hscl\linewidth]{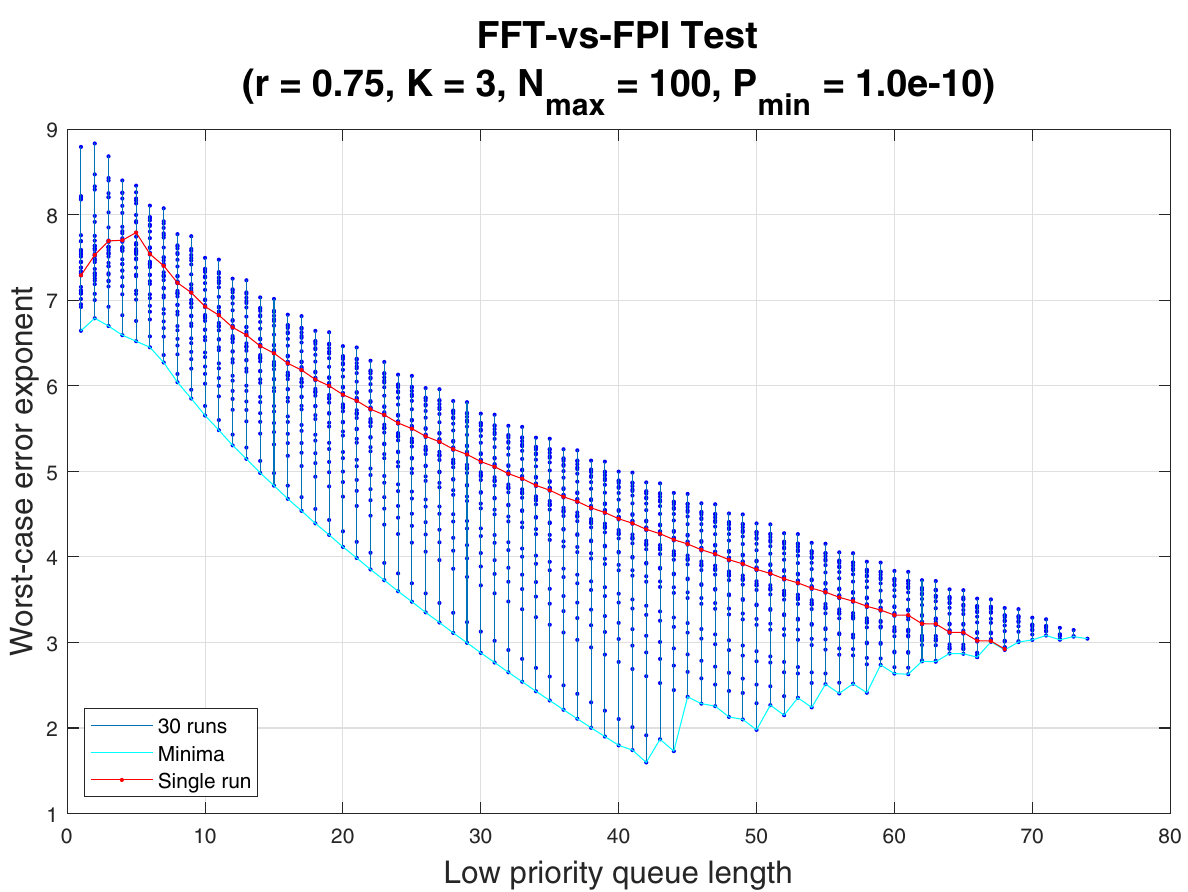}}
{\hphantom{x}\label{TestRepFPIFFTJointLevelPMF}}
{FFT versus FPI test for $K = 3$ priority levels with $30$ randomized level traffic intensities corresponding to
     total traffic intensity $r = 0.75$. The number of decimal places of agreement is
     plotted on the vertical axis as a function of lowest-priority queue length.
     All joint-PMF data points above $P_{\text{min}} = 1.0\times 10^{-10}$ that occur within a
     maximum queue length of $100$ in each priority dimension were considered.}
\end{figure}

\subsection{Results}
\label{Results}
Figures~\ref{TestRepAggFFTJointLevelPMF}--\ref{TestRepNNFFTJointLevelPMF}
present the results of the numerical tests described above.
MOP values prior to worst-case minimization, relevant to the FFT mixture computations,
are displayed on the vertical axes against the relevant queue lengths.
Each test pertains to $K = 4$ priority levels repeated with $30$ randomized level traffic intensities,
all corresponding to total traffic intensity $r = 0.9$, as shown by the blue lines and dots.
The red curve picks out an individual run. The cyan curve delineates the overall worst-case behaviour.
The PMF array is generated up to a maximum queue length of $N_{\text{max}} = 100$
in each priority dimension.
We see that
\mbox{$\Xi_\alpha > 9.5$}
for
$\alpha = \text{agg}, \text{xhi}, \text{xlo},\text{nn}$.

Figure~\ref{TestRepFPIFFTJointLevelPMF} compares the joint PMF arrays as computed
via the FFT and FPI methods, for $K=3$ priority levels with $30$ randomized level
traffic intensities summing to $r = 0.75$. All elements whose probabilities
exceeded a threshold tail level of
\mbox{$P_{\text{min}}> 10^{-10}$}
were included.
The worst case for each lowest-priority queue length is plotted.
We see that \mbox{$\Xi_{\text{fpi}} \gtrsim 2$}. Decreasing performance
the further one progresses into
the tail is a reflection of the truncation error in the FPI method.

\begin{figure}
\FIGURE
{\includegraphics[width=\wscl\linewidth, height=\hscl\linewidth]{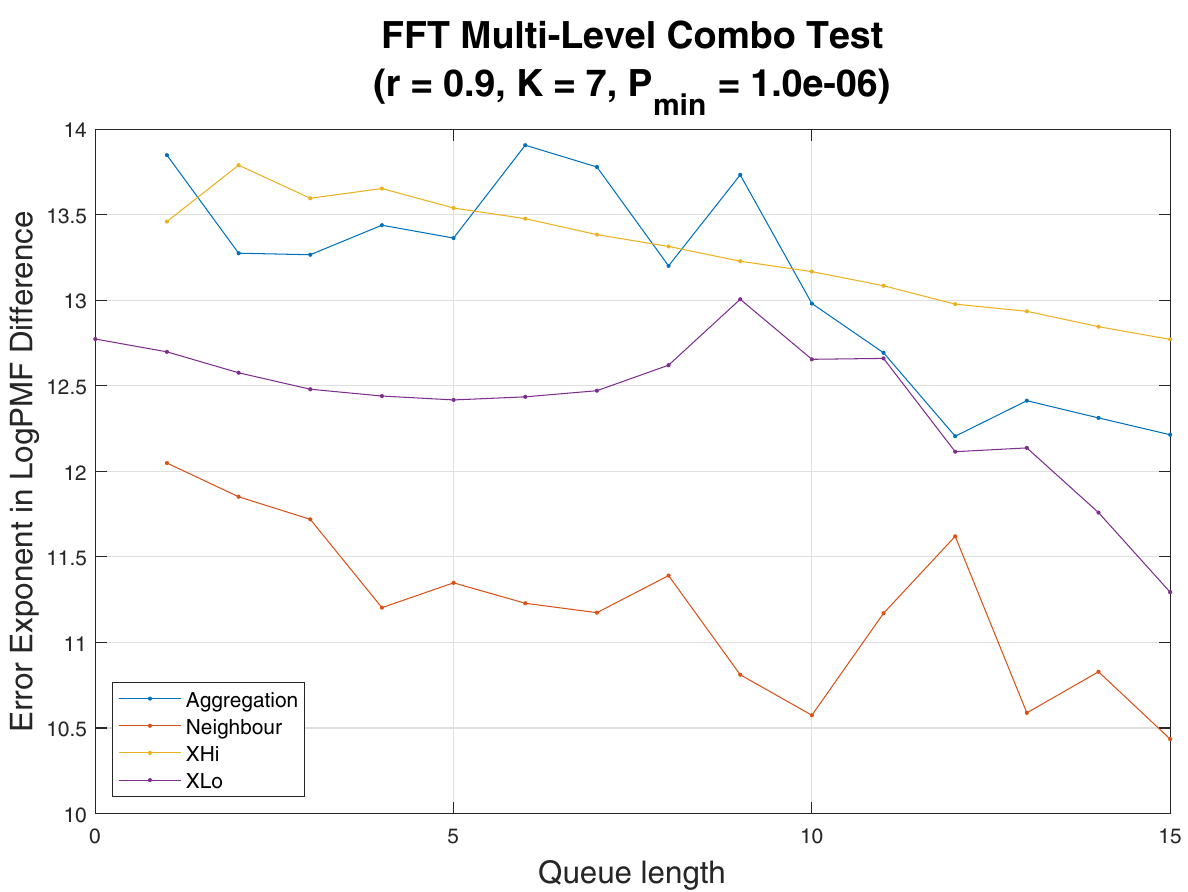}}
{\hphantom{x}\label{TestComboFFTJointLevelPMF}}
{All FFT tests for $K = 7$ priority levels with a random set of level traffic intensities corresponding to
     total traffic intensity $r = 0.9$. The number of decimal places of agreement is
     plotted on the vertical axis as a function of the relevant queue length.
     Joint-PMF data points were computed to a
     maximum queue length of $15$ in each priority dimension.
     For the nearest-neighbour test, all all points above
     $P_{\text{min}} = 1.0\times 10^{-6}$  were considered.}
\end{figure}

Finally, in Figure~\ref{TestComboFFTJointLevelPMF}, we plot the results of the
combined set of tests for $K = 7$ priority levels with a randomly generated
set of level traffic intensities that sum to a total traffic intensity of $r = 0.9$.
The PMF array is generated up to a maximum queue length of $N_{\text{max}} = 15$
in each priority dimension. For the nearest-neighbour test, array elements with
probabilities above
\mbox{$P_{\text{min}} = 10^{-6}$}
are considered, and the results are presented as the worst case for each lowest-priority
queue length. The other curves are self-explanatory.
It can be observed that the algorithm holds up well in this large example, where the
FFT size is 0.25~GB, and the size of the PMF array is 2~GB.

The results obtained confirm (i) numerical stability, (ii) internal consistency,
(iii) consistency between methods, and (iv) agreement with theoretical exact values.

\section{Conclusions}
\label{Concl}
Simple methods for accurate computation of the joint queue-length distribution
for a non-preemptive multi-level priority queue have been developed.
An explicit multi-variate PGF for the joint queue-length distribution
has been derived in closed form for the first time.
With the aid of an FFT mixture method, it is used to numerically compute the
joint queue-length PGF and the queue-length marginals.
A direct method based on a fixed-point iteration of the stationary balance
equations has also been developed. Comparable explicit and exact results
are not available elsewhere in the existing literature.

\begin{APPENDICES}
\section{Multi-Variate PGF}
\label{AppPGF}
In order to solve (\ref{G0U}) for
\begin{equation}
G'_0(u_k,\ldots,u_{K-1}) =  G_0(0,\ldots,0,u_k,\ldots,u_{K-1}) \;,
\end{equation}
we consider a general problem for a family of functions $G_0(.)$ that are differentiated
from one another according to the number of arguments that they take.
The system to be solved has the following structure:
For each
\mbox{$n = 1,2,\ldots,K-1$},
\begin{equation}
G_0(x_1,\ldots,x_n) = \frac{1}{1/x_1 -\lambda_+(x_1,\ldots,x_n)}
     \sum_{k=2}^{n+1}\left(\frac{1}{x_{k-1}} - \frac{1}{x_k}\right)
     G_0(x_k,\ldots,x_n) \;,
\label{G0X}
\end{equation}
with
\mbox{$x_{n+1} \equiv 1$},
and
\mbox{$G_0() = P_0$}.
The functions $\lambda_\pm(x_1,\ldots,x_n)$ are solutions of the quadratic equation
\begin{equation}
\zeta^2 - [1+r -\alpha(x_1,\ldots,x_n)]\zeta + \sigma_{K-n} = 0 \;,
\end{equation}
with
\begin{equation}
\alpha(x_1,\ldots,x_n) = \sum_{k=1}^n x_k r_{K-n+k} \;.
\end{equation}
We observe that
\mbox{$G'_0(u_k,\ldots,u_{K-1}) = G_0(x_1,\ldots,x_n)$},
being the application with
\mbox{$n = K-k$}
arguments,
in which case
\begin{equation}
\alpha(u_k,\ldots,u_{K-1}) = \sum_{\ell=1}^{K-k} u_{\ell+k-1} r_{\ell+k}
     = \sum_{\ell=k}^{K-1} u_\ell r_{\ell+1} \;,
\end{equation}
which involves only the $K-k$ lowest priorities with level traffic intensities
\mbox{$r_{k+1},\ldots,r_K$},
while the constant term of the quadratic
\mbox{$\sigma_{k} = \sum_{\ell=1}^{k}r_\ell$}
aggregates the disjoint highest $k$ priorities with level traffic intensities
\mbox{$r_1, \ldots, r_{k}$}.
This is consistent with the interpretation of
\mbox{$G'_0(u_k,\ldots,u_{K-1})$}
as solving the $(K{+}1{-}k)$-level
sub-problem with the $k$ highest priorities aggregated into a
single level while keeping the total traffic $r$ intensity constant.

Next, we introduce
\begin{align}
\begin{aligned}
J(x_n,\ldots,x_1)         &\equiv G_0(x_1,\ldots,x_n) \;, \\
\zeta_\pm(x_n,\ldots,x_1) &\equiv \lambda_\pm(x_1,\ldots,x_n) \;,
\end{aligned}
\end{align}
and set
\mbox{$z_k \equiv x_{n+1-k}$}
for
\mbox{$k = 1,2,\ldots,n$},
so that
\mbox{$z_1 = x_n$},
\mbox{$z_n = x_1$},
and the convention
\mbox{$x_{n+1} = 1$}
translates to
\mbox{$z_0 = 1$}.
Then (\ref{G0X}) becomes
\begin{align}
\begin{aligned}
J(z_1,\ldots,z_n) &= \frac{1}{1/z_n - \zeta_+(z_1,\ldots,z_n)}
     \sum_{k = 2}^{n+1}\left(\frac{1}{z_{n+2-k}} - \frac{1}{z_{n+1-k}}\right)J(z_1,\ldots,z_{n+1-k}) \\
&= \frac{1}{1/z_n - \zeta_+(z_1,\ldots,z_n)}
     \sum_{\ell = 0}^{n-1}\left(\frac{1}{z_{\ell+1}} - \frac{1}{z_\ell}\right)J(z_1,\ldots,z_\ell) \;,
\end{aligned}
\end{align}
under the change of summation variable
\mbox{$\ell = n + 1 - k$},
and where
\mbox{$J() = G_0(0,\ldots,0) = P_0$}.
Let us now write
\mbox{$J_\ell(\mathbf{z}) \equiv J(z_1,\ldots,z_\ell)$}
and
\mbox{$J_0 \equiv J()$},
so that the subscript $\ell$ indicates the number of elements ({\it i.e.}\ dimension)
of the vector argument $\mathbf{z}$.
Then, we have
\begin{equation}
J_n = \frac{1}{1/z_n - \zeta_+(z_1,\ldots,z_n)}\sum_{\ell=0}^{n-1}
     \left(\frac{1}{z_{\ell+1}} - \frac{1}{z_\ell}\right)J_\ell \;,
\end{equation}
for
\mbox{$n = 1,2,\ldots,K-1$}.
Finally, let
\mbox{$J'_n \equiv \left(1/z_n - \zeta_+(z_1,\ldots,z_n)\right)J_n$}
so that we obtain the recurrence
\begin{equation}
J'_n = \sum_{\ell=0}^{n-1} W_\ell J'_\ell \;,
\end{equation}
with
\begin{equation}
W_\ell = \left(\frac{1}{z_{\ell+1}} - \frac{1}{z_\ell}\right)
     \frac{1}{1/z_\ell - \zeta_+(z_1,\ldots,z_\ell)} \;,
\end{equation}
for
\mbox{$\ell \geq 1$},
and
\mbox{$W_0 = 1/z_1 - 1$}.

For
\mbox{$n \geq 1$},
we can write
\begin{equation}
J'_{n+1} = W_nJ'_n + \sum_{\ell=0}^{n-1} W_\ell J'_\ell = (1 + W_n)J'_n \;,
\end{equation}
which is solved by
\begin{equation}
J'_n = \prod_{\ell=1}^{n-1}(1 + W_\ell){\cdot}J'_1 \;, \quad J'_1 = W_0J'_0 \;,
\end{equation}
where we have set
\mbox{$J'_0 = P_0$}.
After some algebraic manipulation, this leads to the explicit representation
\begin{equation}
J_n(\mathbf{z}) = P_0{\cdot}\prod_{\ell=1}^n\frac{1 - z_\ell\zeta_+(z_1,\ldots,z_{\ell-1})}
     {1 - z_\ell\zeta_+(z_1,\ldots,z_\ell)} \;,
\end{equation}
where the identity
\mbox{$\zeta_+() = 1$}
has been invoked.
The quantity that is ultimately of interest for the $K$-level problem is
\mbox{$J_{K-1}(\mathbf{z}) = G_0(\mathbf{z}) = G_0(z_1,\ldots,z_{K-1})$},
where the arguments of $G_0(.)$ are now reversed.

To complete the specification of the solution,
we recall that
\begin{equation}
\alpha(u_m,\ldots,u_{K-1}) = \sum_{k=m}^{K-1} u_kr_{k+1} \;,
\end{equation}
and set
\mbox{$z_m \equiv u_{K-m}$}
for
\mbox{$m = 1,2\ldots,K-1$},
in order to define the function
\mbox{$\beta(\mathbf{z})$}
via
\begin{equation}
\beta(u_{K-1},\ldots,u_m) \equiv \alpha(u_m,\ldots,u_{K-1}) \;,
\end{equation}
so that
\begin{equation}
\beta(z_1,\ldots,z_{K-m}) = \sum_{k=m}^{K-1} z_{K-k}r_{k+1} \;.
\end{equation}
Equivalently, since $m$ is arbitrary,
\begin{equation}
\beta(z_1,\ldots,z_n) = \sum_{\ell=1}^nz_\ell r_{K+1-\ell} \;.
\end{equation}
The functions $\zeta_\pm(z_1,\ldots,z_n)$ solve the quadratic equation
\begin{equation}
\zeta^2 - [1 + r - \beta(z_1,\ldots,z_n)]\zeta + \sigma_{K-n} = 0  \;.
\end{equation}
The two branches are given by
\begin{equation}
\zeta_\pm(z_1,\ldots,z_n) = \half\left[1+r - \beta(z_1,\ldots,z_n)
     \pm\sqrt{(1+r - \beta(z_1,\ldots,z_n))^2 - 4\sigma_{K-n}}\right] \;,
\end{equation}
for
\mbox{$n = 1,2,\ldots,K-1$}.
Noting that
\mbox{$\sigma_K = \sum_{k=1}^Kr_k = r$},
and with the convention
\mbox{$\beta() \equiv 0$},
we recover the relationships
\mbox{$\zeta_+() = 1$},
\mbox{$\zeta_-() = r$}
from extension to the case
\mbox{$n = 0$}.

\section{Extraction of Marginals}
\label{AppMarg}
In this appendix, we show that all marginals for the general $K$-level priority problem
can be inferred from knowledge of the structure of the PGF for the two-level problem.
It is convenient to introduce the notation
\begin{align}
\begin{aligned}
\xi^\pm_k(\mathbf{z}) &\equiv \zeta_\pm(z_1,\ldots,z_k)\;, \\
\beta_k(\mathbf{z}) &\equiv \beta(z_1,\ldots,z_k)\;,
\end{aligned}
\end{align}
where the integer subscript indicates the number of elements ({\it i.e.}\ the dimension)
of the vector argument.
Then, from (\ref{PGFzp}) and (\ref{GMinus}), it follows that
\begin{equation}
\sum_{\ell=0}^\infty G_\ell(\mathbf{z}) = \sum_{k=0}^{K-1}\frac{1 - z_k\zeta^-_k(\mathbf{z})}
     {1 - z_{k+1}\zeta^-_k(\mathbf{z})} \;,
\label{GSum}
\end{equation}
where we have set
\mbox{$z_K \equiv 1$},
and noted that
\begin{equation}
P_0 = 1 - r = 1 - \zeta_-() = 1 - \zeta^-_0(\mathbf{z}) \;.
\end{equation}
To extract the $p$-th marginal, we must consider the vector argument
\mbox{$\mathbf{z} = \mathbf{z}^{(p)}$}
where we set
\mbox{$z_k = 1$}
for all
\mbox{$k \neq p$}.
Thus
\mbox{$\mathbf{z}^{(p)} \equiv \mathbf{1} + (z_p-1)\mathbf{e}_p$}.
We may observe that the $k$-th term of the product on the RHS of (\ref{GSum})
is equal to unity, expect when
\mbox{$k = p-1$}
or
\mbox{$k = p$}.
This leads to a representation of the PGF for the $p$-th marginal as
\begin{equation}
G_{\text{mrg}}^{(p)}(z_p) = \sum_{\ell=0}^\infty G_\ell(\mathbf{z}^{(p)}) =
     \frac{1 - \zeta^-_{p-1}(\mathbf{z}^{(p)})}{1 - z_p\zeta^-_{p-1}(\mathbf{z}^{(p)})}{\cdot}
     \frac{1 - z_p\zeta^-_{p}(\mathbf{z}^{(p)})}{1 - \zeta^-_{p}(\mathbf{z}^{(p)})} \;,
\end{equation}
for
\mbox{$p = 1,2,\ldots,K-1$}.
Since the identity
\mbox{$\beta_k(\mathbf{1}) = r - \sigma_{K-k}$}
implies that
\mbox{$\zeta^+_k(\mathbf{1}) = 1$},
\mbox{$\zeta^-_k(\mathbf{1}) = \sigma_{K-k}$},
we have
\begin{equation}
\zeta^-_{p-1}(\mathbf{z}^{(p)}) = \zeta^-_{p-1}(\mathbf{1}) = \sigma_{K+1-p} \;.
\end{equation}
Next, we observe that
\begin{equation}
\beta_p(\mathbf{z}^{(p)}) = \beta_p(1,\ldots,z_p,\ldots,1)
     = \sum_{k=1}^{p-1}r_{K+1-k} + r_{K+1-p}z_p \;.
\end{equation}
Setting
\begin{equation}
r_{\text{lo}} = r_{K+1-p} \;, \quad r_{\text{hi}} = \sum_{k=1}^{K-p} r_k  = \sigma_{K-p}\;,
     \quad r_{\text{sum}} \equiv r_{\text{lo}} + r_{\text{hi}} = \sigma_{K+1-p}\;,
\end{equation}
we obtain
\begin{equation}
r - \beta_p(\mathbf{z}^{(p)}) = r_{\text{sum}} - r_{\text{lo}} z_p \;.
\end{equation}
Therefore,
\begin{align}
\begin{aligned}
\zeta^\pm_p(\mathbf{z}^{(p)}) &= \half\left[1 + r - \beta_p(\mathbf{z}^{(p)})
     \pm \sqrt{\left(1 + r - \beta_p(\mathbf{z}^{(p)})\right)^2 - 4\sigma_{K-p}}\right] \\
&=   \half\left[1 + r_{\text{sum}} - r_{\text{lo}}z_p
     \pm \sqrt{(1 + r_{\text{sum}} - r_{\text{lo}}z_p)^2 - 4r_{\text{hi}}}\right] \\
&= \zeta_\pm(z_p) \;,
\end{aligned}
\end{align}
where the final identification with
\mbox{$\zeta_\pm(z)$}
for the two-level problem is associated with the mappings
\mbox{$r_1 \leftarrow r_{\text{hi}}$},
\mbox{$r_2 \leftarrow r_{\text{lo}}$},
\mbox{$r \leftarrow r_{\text{sum}}$}.
Consequently,
\begin{equation}
G_{\text{mrg}}^{(p)}(z_p) = \frac{1 - r_{\text{sum}}}{1 - r_{\text{sum}}z_p}{\cdot}
     \frac{1 - z_p\zeta_-(z_p)}{1 - \zeta_-(z_p)} \;,
\end{equation}
for
\mbox{$p = 1,2,\ldots,K-1$}.
Hence, on comparing with (\ref{GMinus}) for the case
\mbox{$K = 2$},
we see that the desired result for the marginal PGFs is established.
Alternatively, by multiplying numerator and denominator by
\mbox{$1 - \zeta_+(z_p)$}
and using the identities given in (\ref{ZetaIds}), one recovers the first form in (\ref{GLo2}).

\section{PGF Ratio}
\label{AppRatio}
Considering the quantities
\mbox{$\mathscr{P}^\pm_\kappa(\mathbf{z})$},
\mbox{$\mathscr{Q}^\pm_\kappa(\mathbf{z})$}
as defined in (\ref{PQ}),
we shall show that
\begin{equation}
\mathscr{P}^+_\kappa(\mathbf{z})\mathscr{P}^-_\kappa(\mathbf{z}) =
     \mathscr{Q}^+_\kappa(\mathbf{z})\mathscr{Q}^-_\kappa(\mathbf{z}) \;,
\end{equation}
which implies that one can eliminate the removable singularities in (\ref{G0PQ})
by making the substitution
\begin{equation}
\mathscr{P}^+_\kappa(\mathbf{z})/\mathscr{Q}^+_\kappa(\mathbf{z}) \mapsto
     \mathscr{Q}^-_\kappa(\mathbf{z})/\mathscr{P}^-_\kappa(\mathbf{z}) \;.
\end{equation}
In order to first evaluate
\mbox{$\mathscr{P}^+_\kappa\mathscr{P}^-_\kappa$},
we let
\mbox{$a \equiv 4\sigma_{K+1-\kappa}$}
and
\mbox{$b \equiv 1 + r - \beta_{\kappa-1}$}.
Then it is easy to see that
\begin{align}
\begin{aligned}
\mathscr{P}^+_\kappa\mathscr{P}^-_\kappa &= 1 - bz_\kappa + \frac{a}{4}z^2_\kappa \\
&= 1 - (1 + r -\beta_{\kappa-1})z_\kappa + \sigma_{K+1-\kappa}z^2_\kappa \;.
\end{aligned}
\end{align}
Similarly, if we now let
\mbox{$a' \equiv 4\sigma_{K-\kappa}$}
and
\mbox{$b' \equiv 1 + r - \beta_\kappa$},
then
\begin{align}
\begin{aligned}
\mathscr{Q}^+_\kappa\mathscr{Q}^-_\kappa &= 1 - b'z_\kappa + \frac{a'}{4}z^2_\kappa \\
&= 1 - (1 + r -\beta_\kappa)z_\kappa + \sigma_{K-\kappa}z^2_\kappa \\
&= 1 - (1 + r -\beta_{\kappa-1} - r_{K+1-\kappa}z_\kappa)z_\kappa + \sigma_{K-\kappa}z^2_\kappa \\
&= 1 - (1 + r -\beta_{\kappa-1})z_\kappa + (\sigma_{K-\kappa} + r_{K+1-\kappa})z^2_\kappa \\
&= 1 - (1 + r -\beta_{\kappa-1})z_\kappa + \sigma_{K+1-\kappa}z^2_\kappa \;.
\end{aligned}
\end{align}
The desired result follows from the equality of the right-hand sides of the foregoing pair of equations.

\section{Alternative Derivation of the Joint PGF}
\label{AltMulti}
Recalling from (\ref{FullPGF}) that
\mbox{$G_\ell(z_1,\ldots,z_{K-1}) = G_0(z_1,\ldots,z_{K-1})\zeta_-^\ell$},
we obtain for the full $K$-dimensional joint PGF,
\begin{equation}
G(z_1,\ldots,z_K) \equiv \sum_{\ell=0}^\infty z_K^\ell G_\ell(z_1,\ldots,z_{K-1})
     = \frac{G(z_1,\ldots,z_{K-1},0)}{1 - z_K\zeta_{K-1}(\mathbf{z})} \;,
\label{GK}
\end{equation}
where
\mbox{$\zeta_{K-1}(\mathbf{z}) \equiv \zeta_-(z_1,\ldots,z_{K-1})$}
such that the subscript denotes the number of arguments.
We introduce the notation
\mbox{$G(z_1,\ldots,z_p|z) \equiv G(z_1,\ldots,z_p,z,\ldots,z)$}
for any
\mbox{$z\in\mathbb{C}$}
and
\mbox{$p = 0,1,\ldots,K$}.
For
\mbox{$p = 0$},
we have
\mbox{$G({\cdot}|z) = G(z,\ldots,z)$},
while for
\mbox{$p = K$},
we have
\mbox{$G(z_1,\ldots,z_K|z) = G(z_1,\ldots,z_K)$},
independent of $z$.

Since the marginal distribution of the aggregation of the top $p$ priority levels is also geometric,
as argued in the discussion below (\ref{Gaggp}),
we have that
\begin{equation}
G(z_1,\ldots,z_{K-p}|z) = \frac{G(z_1\ldots,z_{K-p}|0)}{1 - z\zeta_{K-p}(\mathbf{z})} \;,
\label{GKp}
\end{equation}
for arbitrary $z$
and
\mbox{$p = 1,2,\ldots,K$},
where the rates $\zeta_{K-p}(\mathbf{z})$
are to be determined.
When
\mbox{$p = 1$},
we recover (\ref{GK}).
Setting
\mbox{$z = z_{K-p+1}$} in (\ref{GKp})
yields
\begin{equation}
G(z_1,\ldots,z_{K-p}|z_{K-p+1}) = \frac{G(z_1\ldots,z_{K-p}|0)}{1 - z_{K-p+1}\zeta_{K-p}(\mathbf{z})} \;,
\label{GKp1}
\end{equation}
for
\mbox{$p = 1,\ldots,K$}.
Next we write (\ref{GKp}) as
\begin{equation}
G(z_1,\ldots,z_{K-p+1}|z) = \frac{G(z_1\ldots,z_{K-p+1}|0)}{1 - z\zeta_{K-p+1}(\mathbf{z})} \;,
\label{GKp2}
\end{equation}
and choose
\mbox{$z_{K-p+1} = z$}
so that
\begin{equation}
G(z_1,\ldots,z_{K-p}|z) = \frac{G(z_1\ldots,z_{K-p+1}|0)}{1 - z\zeta_{K-p+1}(z_1,\ldots,z_{K-p},z)} \;.
\label{GKp3}
\end{equation}
Now we set
\mbox{$z = z_{K-p+1}$}
to obtain
\begin{equation}
G(z_1,\ldots,z_{K-p}|z_{K-p+1}) = \frac{G(z_1\ldots,z_{K-p+1}|0)}{1 - z_{K-p+1}\zeta_{K-p+1}(\mathbf{z})} \;.
\label{GKp4}
\end{equation}
We observe that (\ref{GKp1}) and (\ref{GKp4}) express equations for the same quantity.
Thus, on equating the RHS of each equation,
we arrive at the recursion
\begin{equation}
G(z_1,\ldots,z_k|0) = \frac{1 - z_k\zeta_k(\mathbf{z})}{1 - z_k\zeta_{k-1}(\mathbf{z})}{\cdot}
     G(z_1,\ldots,z_{k-1}|0) \;,
\end{equation}
for
\mbox{$k = 1,2,\ldots,K-1$},
which leads to the result
\begin{equation}
G_0(z_1,\ldots,z_{K-1}) = G(z_1,\ldots,z_{K-1}|0)
     = P_0\prod_{k=1}^{K-1}\frac{1 - z_k\zeta_k(\mathbf{z})}{1 - z_k\zeta_{k-1}(\mathbf{z})} \;,
\end{equation}
where we have used the fact that
\mbox{$P_0 = G({\cdot}|0) = G(0,\ldots,0)$}.

In order to determine the rates
$\zeta_k(\mathbf{z})$,
we multiply the stationary balance equation (\ref{Vector}) by
\mbox{$z_1^{n_K}\cdots z_{K-k}^{n_{k+1}}{\cdot}z^{n_k+\cdots+n_1}$},
and sum over all indices
\mbox{$n_1,\ldots n_K$}
to obtain
\begin{equation}
\left[\beta_k(\mathbf{z}) - r +\frac{\sigma_{K-k}}{z\zeta_k(\mathbf{z})}\right]G(z_1,\ldots,z_k|z)
     = G(z_1,\ldots,z_k|0) \;,
\end{equation}
where
\begin{equation}
\beta_k(\mathbf{z}) = \sum_{\kappa=1}^k z_\kappa r_{K-k+1} \;, \quad
  \sigma_k = \sum_{\kappa = 1}^k r_\kappa \;,
\end{equation}
and we have used the relation
\begin{equation}
\left.\frac{d}{dz}G(z_1,\ldots,z_k|z)\right|_{z=0} = \zeta_k(\mathbf{z}) G(z_1,\ldots,z_k|0) \;,
\end{equation}
that follows from (\ref{GKp}).
We also use (\ref{GKp}), expressed as
\begin{equation}
G(z_1,\ldots,z_k|0) = (1 - z\zeta_k(\mathbf{z}))G(z_1,\ldots,z_k|z) \;,
\end{equation}
to obtain
\begin{equation}
\left[z^2\zeta^2_k - (1 + r - \beta_k)z\zeta_k + \sigma_{K-k}\right]
     G(z_1,\ldots,z_k|z) = 0 \;,
\end{equation}
provided
\mbox{$z\neq 0$}.
Setting
\mbox{$z = 1$}
yields
\begin{equation}
\zeta^2_k - (1 + r - \beta_k)\zeta_k + \sigma_{K-k} = 0 \;,
\end{equation}
which is solved by (\ref{ZetaBeta}).
Therefore, the result (\ref{GMinus}) is established.
\end{APPENDICES}


%
\section*{Acknowledgments}
The authors gratefully acknowledge useful discussions with Dr.~Stephen Bocquet.


\bibliographystyle{informs2014} 
\bibliography{MultiNP} 



\end{document}